\documentclass[10pt, a4paper]{amsart}
\usepackage{amsfonts,amsmath,latexsym,amssymb}

\def\C{\mathbb{C}}
\def\N{\mathbb{N}}
\def\P{\mathbb{P}}%
\def\R{\mathbb{R}}
\def\Z{\mathbb{Z}}
\def\Bs{\mathop{\mathrm{Bs}}\nolimits}
\def\Cls{\mathop{\mathrm{Cls}}\nolimits}
\def\id{\mathop{\mathrm{id}}\nolimits}
\def\Int{\mathop{\mathrm{Int}}\nolimits}
\def\Item#1{\par\hangindent\parindent\indent\llap{#1\enspace}\ignorespaces}
\def\Ker{\mathop{\mathrm{Ker}}\nolimits}
\def\Krulldim{\mathop{\rm Krull \; dim}\nolimits}
\def\length{\mathop{\mathrm{length}}\nolimits}
\def\lto{\longrightarrow}
\def\mult{\mathop{\mathrm{mult}}\nolimits}

\def\Rat{\mathop{\rm Rat}\nolimits}
\def\Reg{\mathop{\mathrm{Reg}}\nolimits}
\def\sign{\mathop{\mathrm{sign}}\nolimits}
\def\Sing{\mathop{\mathrm{Sing}}\nolimits}
\def\subsetne{\subsetneq}
\def\supsetne{\supsetneq}
\def\XYZ{\mathop{\mathcal{X}}\nolimits}
\def\Zar{\mathop{\mathrm{Zar}}\nolimits}
\def\cC{\mathcal{C}}
\def\cE{\mathcal{E}}
\def\cH{\mathcal{H}}
\def\cI{\mathcal{I}}
\def\cJ{\mathcal{J}}
\def\cL{\mathcal{L}}
\def\cP{\mathcal{P}}
\def\cQ{\mathcal{Q}}
\def\cR{\mathcal{R}}
\def\cZ{\mathcal{Z}}
\def\frf{\mathfrak{f}}
\def\frg{\mathfrak{g}}
\def\frh{\mathfrak{h}}

\catcode`\@=11
\outer\def\proclaim#1{%
  \removelastskip\penalty-400\vskip0.8em plus0.3em minus0.3em
  {\bf#1.}}%
\def\endproclaim{\par\penalty-400\vskip0.8em plus0.3em minus0.3em}%
\def\Proof{\removelastskip\par\begin{proof}}%

\catcode`\@=11
\newdimen\htqed \newdimen\wdqed \newdimen\dpqed
\def\hidehrule#1#2{\kern-#1 \hrule height#1 depth#2 \kern-#2 }
\def\hidevrule#1#2{\kern-#1{\dimen0=#1
    \advance\dimen0 by#2\vrule width\dimen0}\kern-#2 }
\def\makeblankbox#1#2{\hbox{\lower\dpqed\vbox{\hidehrule{#1}{#2}\kern-#1 %
    \hbox to \wdqed{\hidevrule{#1}{#2}%
    \raise\htqed\vbox to #1{}%
    \lower\dpqed\vtop to #1{}%
    \hfil\hidevrule{#2}{#1}}%
    \kern-#1\hidehrule{#2}{#1}}}}
\def\QED{\htqed=6.7pt \dpqed=0pt \wdqed=6.7pt
    {\unskip\nobreak\hfil\penalty50\quad\null\nobreak\hfil
    {\hbox{\makeblankbox{0.17pt}{0.17pt}}}
    \parfillskip0pt\finalhyphendemerits0\par\medskip}}

\catcode`\@=\active
\catcode`\@=11

\def\Ttac{Theorem 1.1}
\def\Ttad{Theorem 1.2}
\def\Ttadn{1.2}
\def\Ttadb{Theorem 1.4}
\def\Ttadbn{1.4}
\def\Ttae{Theorem 1.5}
\def\Ttaf{Theorem 1.6}
\def\Ttab{Theorem 1.7}
\def\Tqbv{Definition 2.2} 
\def\Tqba{Definition 2.6} 
\def\Tqbw{Proposition 2.7} 
\def\Tqbb{Proposition 2.8} 
\def\Tqbc{Definition 2.9} 
\def\Tqbf{Lemma 2.10} 
\def\Tqbae{Theorem 2.11} 
\def\Tqbg{Definition 2.12} 
\def\Tqbh{Proposition 2.13} 
\def\Tqbi{Example 2.14} 
\def\Tqbaib{Remark 2.15} 
\def\Tqacnode{Definition 2.16} 
\def\Tqbaf{Theorem 2.17} 
\def\Tqbag{Theorem 2.18} 
\def\Ttaaa{Proposition 3.1}
\def\Ttaca{Proposition 3.3}
\def\Ttada{Proposition 3.4}
\def\Ttadan{3.4}
\def\Ttaea{Lemma 3.5}
\def\Ttaean{3.5}
\def\Ttafb{Lemma 3.6}
\def\Ttafbn{3.6}
\def\Ttagb{Lemma 3.7}
\def\Ttag{Lemma 3.8}
\def\Ttagc{Notation 3.9}
\def\Ttaic{Lemma 3.10}
\def\Ttaicn{3.10}
\def\Ttah{Lemma 3.11}
\def\Ttahn{3.11}
\def\Ttca{Lemma 3.12}
\def\Ttcb{Theorem 3.13}
\def\Ttcbn{3.13}
\def\Ttakb{Lemma 3.14}
\def\Ttaka{Theorem 3.15}
\def\Ttakan{3.15}
\def\Ttala{Definition 3.16}
\def\Ttalb{Lemma 3.17}
\def\Ttalc{Lemma 3.18}
\def\Ttamb{Lemma 3.19}
\def\Ttama{Theorem 3.20}
\def\Ttaman{3.20}
\def\Ttapa{Corollary 3.21}
\def\Ttapb{Theorem 3.22}
\def\Ttapc{Theorem 3.23}
\def\Ttapd{Theorem 3.24}
\def\Ttanc{Corollary 3.25}
\def\Ttana{Proposition 3.26}
\def\Ttanb{Corollary 3.27}
\def\Ttda{Theorem 4.1}

\begin{document}
\title[Extremal Cubic Inequalities]{Extremal Cubic Inequalities of Three Variables}

\vskip2mm

\author{Tetsuya ANDO}

\address{
Department of Mathematics and Informatics, 
Chiba University, \\
Yayoi-cho 1-33, Inage-ku, 
Chiba 263-8522, JAPAN \\}
\email{ando@math.s.chiba-u.ac.jp}

\date{04.03.2022}  



\begin{abstract}
Let $\mathcal{H}_{3,d}$ be the vector space of homogeneous three variable 
polynomials of degree $d$, 
and $\mathcal{P}_{3,d}^+$ be the set of all elements $f \in \mathcal{H}_{3,d}$ 
such that $f(x,y,z) \geq 0$ for all $x \geq 0$, $y \geq 0$, $z \geq 0$. 
In this article we determine all extremal elements of $\mathcal{P}_{3,3}^+$. 
We prove that if $f \in \mathcal{P}_{3,3}^+$ is an irreducible extremal 
element, 
then the zero locus $V_{\mathbb{C}}(f)$ in $\mathbb{P}_{\mathbb{C}}^2$ is 
a rational curve whose unique singularity is an acnode in the interior 
of $\mathbb{P}_+^2$ or a cusp on an edge of $\mathbb{P}_+^2$. 
We also prove that if $f \in \mathcal{P}_{3,3}^+$ is an extremal element, 
then $f(x^2,y^2,z^2)$ is an extremal element of $\mathcal{P}_{3,6}$, 
where $\mathcal{P}_{3,d}$ is the set of all the 
elements $f \in \mathcal{H}_{3,d}$ such that $f(x,y,z) \geq 0$ for all $x$, 
$y$, $z \in \mathbb{R}$. 
A notion of infinitely near zeros of an inequality is introduced, 
and plays an important role. 
\end{abstract}

\maketitle


\section{Introduction}

In this article, we determine all the extremal elements of the convex cone 
$\cP_{3,3}^+$, where basic symbols are defined as the following: 
\begin{align*}
 \R_+ 
 & := \big\{ a \in \R \; \big| \; a \geq 0 \big\}, \\
 \P_+^n 
 & := \big\{(x_0\colon\cdots\colon x_n) \in \P_{\R}^n \; \big| \; 
      \hbox{$x_i x_j \geq 0$ for all $0 \leq i<j \leq m$} \big\}, \\
 \cH_{n,d} 
 & := \R[x_1,\ldots,x_n]_{(d)} \\
 & = \big\{ f \in \R[x_1,\ldots,x_n] \; \big| \; \hbox{$f$ is 
      homogeneous and $\deg f = d$}\big\} \cup \{0\}, \\
 \cP_{n,d} 
 & := \big\{ f \in \cH_{n,d} \; \big| \; \hbox{
     $f(x_1,\ldots,x_n) \geq 0$ for all 
      $(x_1,\ldots, x_n) \in \R^n$}\big\}, \\
 \cP_{n,d}^+ 
 & := \big\{ f \in \cH_{n,d} \; \big| \; \hbox{
     $f(x_1,\ldots,x_n) \geq 0$ for all 
      $(x_1,\ldots, x_n) \in \R_+^n$}\big\}, \\
 \Sigma_{n,d} 
 & := \big\{ f \in \cP_{n,d} \; \big| \; \hbox{
     $f$ is a sum of squares of some elements from $\cH_{n,d/2}$}\big\}. 
\end{align*}
Let $\cP$ be a closed convex cone which contains no lines, 
and $\cH$ be the vector space spanned by $\cP$. 
An element $0 \ne f \in \cP$ is said to be {\it extremal} 
if $g$, $h \in \cP$ and $f=g+h$, 
implies that $g$ and $h$ are divisible by $f$. 
Let 
\[\cE(\cP) := \big\{ f \in \cP \; \big| \; 
  \hbox{$f$ is an extremal element of $\cP$} \big\}.\]
An element $f \in \cE(\cP)$ is said to be {\it exposed} if there exists a 
hyperplane $H \subset \cH$ such that $H \cap \cP = \R_+ \cdot f$. 
For $f \in \cH_{n,d}$ and $K = \R$ or $\C$, we denote 
\[V_K(f) := \big\{ (x_1 \colon \cdots \colon x_{n}) \in \P_K^{n-1} \; 
  \big| \; f(x_1,\ldots,x_n) = 0 \big\}.\]
A set $V_{\R}(f) \cap \P_+^{n-1}$ is denoted by $V_+(f)$. 
We will see that the extremal ray $\R_+ \cdot f$ of $\cP_{3,3}^+$ is 
determined by $V_+(f)$ with additional information of infinitely near points. 

Before we state our main result, i.e., \Ttae, 
we present some elements of $\cE(\cP_{3,3}^+)$. 
For example, $\frf_0(x$, $y$, $z) = x^2y+y^2z+z^2x-3xyz 
  \in \cE(\cP_{3,3}^+)$ (see \Ttana \ 
or \cite[Corollary 3.3]{RefCL}). 
The Schur's inequality type polynomial $\frf_1(x$, $y$, $z) 
 = x^3+y^3+z^3 + 3xyz - x^2y - y^2z - z^2x - xy^2 - yz^2 - zx^2$ 
and $xyz$ are also elements of $\cE(\cP_{3,3}^+)$ (see \Ttanb). 
If $f \in \cE(\cP_{3,3}^+)$ is a symmetric polynomial, 
then $f = \alpha \frf_1$ or $f = \alpha xyz$ ($\exists \alpha>0$. 
See also \cite[Theorem 3.8]{RefCL}). 
On the other hand, AM-GM inequality type polynomial 
$$x^3+y^3+z^3-3xyz = \frf_1(x, y, z) 
 + \frf_0(x, y, z) + \frf_0(y, x, z)$$
is not extremal in $\cP_{3,3}^+$. 

\Ttac, \Ttadn, \Ttadbn \ below describe three new families of 
elements of $\cE(\cP_{3,3}^+)$. 

\def\Ttac{Theorem 1.1}
\proclaim{Theorem 1.1} 
{\sl Assume that $p \geq 0$, $q \geq 0$, $r \geq 0$, $pq-p+1 > 0$, 
$qr-q+1 > 0$ and $rp-r+1 > 0$. Then: } 
{\parindent=20pt
\Item{\rm (1)} There exists an irreducible polynomial 
$\frf_{pqr} \in \cE(\cP_{3,3}^+)$ which satisfies 
\[\frf_{pqr}(1,1,1) = \frf_{pqr}(0,p,1) 
= \frf_{pqr}(1,0,q) = \frf_{pqr}(r,1,0) = 0.\]
\Item{(2)} If $p>0$, $q>0$, $r>0$, 
then every $f \in \cP_{3,3}^+$ satisfying 
$f(1,1,1) = f(0,p,1) = f(1,0,q) = f(r,1,0) = 0$ 
is of the form $f = \alpha \frf_{pqr}$ for some $\alpha \in \R$.} 
\endproclaim

For the explicit expression of $\frf_{pqr}$ see \Ttala. 

\def\Ttad{Theorem 1.2}
\def\Ttadn{1.2}
\proclaim{Theorem 1.2} 
{\sl For $p$, $q \in \R$, let 
\begin{align*}
 & \frg_{pq}(x,y,z) := z^3 + q^2 x^2z + p^2 y^2z - 2q x z^2 - 2p y z^2 
             - (p^2+q^2-4p-4q+3)x y z \\
 & \hskip70pt + (1-p+q)(1-p-q)x^2y + (1+p-q)(1-p-q) x y^2.
\end{align*}
If $p \geq 0$, $q \geq 0$ and $p+q \leq 1$, 
then $\frg_{pq} \in \cE(\cP_{3,3}^+)$ and 
\[\frg_{pq}(1,1,1) = \frg_{pq}(0,1,p) = \frg_{pq}(1,0,q) 
= \frg_{pq}(1,0,0) = \frg_{pq}(0,1,0) = 0.\]
Conversely, if $p>0$, $q>0$, $p+q<1$ and $f \in \cP_{3,3}^+$ satisfies 
$f(1,1,1) = f(1,0,p) = f(0,1,q) = f(1,0,0) = f(0,1,0) = 0$, 
then there exists $\alpha \in \R$ such that $f = \alpha \frg_{pq}$. } 
\endproclaim

\proclaim{Remark 1.3} 
{\rm $\frf_{ppp}$ was discovered in \cite[Theorem 3.1]{RefAa}. 
A special type of $\frg_{pq}$ was discovered in \cite{RefR}. 
Let 
\begin{align*}
 M_t(x,y,z) 
  & := (1-2t^2)(x^4y^2+x^2y^4) + t^4(x^4z^2+y^4z^2) \\
  & \hskip30pt - (3-8t^2+2t^4)x^2y^2z^2 
     - 2t^2(x^2+y^2)z^4 + z^6
\end{align*}
be the polynomial of (1.8) or (6.17) in \cite{RefR}. 
Then $M_t(x$, $y$, $z) = \frg_{t^2,t^2}(x^2$, $y^2$, $z^2)$. }
\endproclaim

Note that $V_{\C}(\frf_{pqr})$ and $V_{\C}(\frg_{pq})$ have a node 
at $(1 \colon 1 \colon 1) \in \P_{\C}^2$. 
There are extremal elements $f \in \cE(\cP_{3,3}^+)$ such 
that $V_{\C}(f)$ has a cusp. 

\def\Ttadb{Theorem 1.4}
\def\Ttadbn{1.4}
\proclaim{Theorem 1.4} 
{\sl Let 
\begin{align*}
 \frh_{pq}(x,y,z) 
 & := 2x^3 + 3(p-q) x^2 y
          - 6p q x y^2 + q^2(3p+q)y^3 
    + 3(q-p)x^2z - 6p q x z^2 \\
 & \hskip20pt 
    + (p^3+36p^2q-6p q^2-2q^3)y^2 z 
    + (-2p^3-6p^2q+3p q^2+q^3)y z^2 \\
 & \hskip20pt 
    + p^2(p+3q) z^3 + 12p q x y z. 
\end{align*}
Assume that $p \geq 0$, $q \geq 0$ and $(p$, $q) \ne (0$, $0)$. 
Then $\frh_{pq} \in \cE(\cP_{3,3}^+)$ and 
$\frh_{pq}(p$, $0$, $1) = \frh_{pq}(q$, $1$, $0) = 0$. 
Moreover, $V_{\C}(f)$ has a cusp at $(0 \colon 1 \colon 1)$. 
Conversely, for $f \in \cP_{3,3}^+$, 
if $V_{\C}(f)$ has a cusp at $(0 \colon 1 \colon 1)$ and 
$f(p,0,1) = f(q,1,0) = 0$, 
then there exists $\alpha \in \R$ such that $f = \alpha \frh_{pq}$. } 
\endproclaim

Our main theorem is the characterization of the elements from 
$\cE(\cP_{3,3}^+)$: 

\def\Ttae{Theorem 1.5}
\proclaim{Theorem 1.5} 
{\sl Let $f(x,y,z) \in \cE(\cP_{3,3}^+)$. 
Then, $f(x,y,z)$ is a positive multiple of one of the following 
polynomials: }\par
{\parindent=20pt
\Item{(1)} {\sl $\frf_{pqr}(x$, $y/a$, $z/b)$ where $a>0$, $b>0$, 
$p \geq 0$, $q \geq 0$, $r \geq 0$, $p q-p+1>0$, $q r-q+1>0$ 
and $r p-r+1>0$.}
\Item{(2)} {\sl $\frf_{prq}(x$, $z/b$, $y/a)$ where $a>0$, $b>0$, 
$p \geq 0$, $q \geq 0$, $r \geq 0$, $p q-q+1>0$, $q r-r+1>0$, 
$r p-p+1>0$ and $pqr=0$.}
\Item{(3)} {\sl $\frg_{pq}(x$, $y/a$, $z/b)$ or 
$\frg_{pq}(y/a$, $z/b$, $x)$ or $\frg_{pq}(z/b$, $x$, $y/a)$ 
where $a>0$, $b>0$, $p \geq 0$, $q \geq 0$ and $p+q<1$. }
\Item{(4)} {\sl $\frh_{pq}(x$, $y/a$, $z)$ or 
$\frh_{pq}(y$, $z/a$, $x)$ 
or $\frh_{pq}(z$, $x/a$, $y)$ where $a>0$, $p \geq 0$, $q \geq 0$ 
and $(p$, $q) \ne (0$, $0)$. }
\Item{(5)} {\sl $x(ax+by+cz)^2$ or $y(ax+by+cz)^2$ or $z(ax+by+cz)^2$ 
where $a$, $b$, $c \in \R$, $(a$, $b$, $c) \ne (0$, $0$, $0)$ and 
$\dim \big(V_{\R}(ax+by+cz) \cap \P_+^2\big) = 1$.}
\Item{(6)} {\sl The monomial $xyz$.}

}
{\sl Conversely, all polynomials in (1)-(6) belongs to $\cE(\cP_{3,3}^+)$. }
\endproclaim

This theorem will be proved in \S 3.6. 

Hilbert proved that $\cP_{n,d} = \Sigma_{n,d}$ if and only if 
$n \leq 2$ or $d=2$ or $(n$, $d)=(3$, $4)$. 
Moreover, every element of $\cE(\cP_{3,4})$ is 
a square of a quadratic polynomial \cite{RefHil}. 
We shall give an alternative proof of this fact at \Ttda. 

The first part of the following theorem is proved 
in \cite[Remark 8]{RefBHORS}, 
and the second part follows from \cite[Theorem 7.2]{RefR}. 

\def\Ttaf{Theorem 1.6}
\proclaim{Theorem 1.6} 
{\sl If $f \in \cP_{3,6}$ is an exposed extremal element which is not 
the square of a cubic polynomial, then $V_{\C}(f)$ is an irreducible 
rational curve which has 10 acnodes $P_1$,$\ldots$, $P_{10}$, and 
$V_{\R}(f) = \{P_1$,$\ldots$, $P_{10}\}$. 
On the other hand, 
if $f \in \cP_{3,6}$ and $V_{\C}(f)$ is an irreducible curve 
which has 10 acnodes in $\P_{\R}^2$, 
then $f \in \cE(\cP_{3,6})$. }
\endproclaim

In spite of this general theorem, only a few concrete elements of 
$\cE(\cP_{3,6})$ ware known (see \cite{RefKS}). 
But as a corollary of \Ttae, we obtain the following: 

\def\Ttab{Theorem 1.7}
\proclaim{Theorem 1.7} 
{\sl If $f(x,y,z) \in \cE(\cP_{3,3}^+)$, 
then $f(x^2,y^2,z^2) \in \cE(\cP_{3,6})$. }
\endproclaim

The convex cone $\cP_{3,6}$ is studied 
in \cite{RefCL, RefCLR, RefCLRb, RefCLRc, RefKS, RefR}, 
and the convex cone $\cP_{3,6}^e
  := \cP_{3,6} \cap \R[x^2,y^2,z^2]$ is studied in \cite{RefGKR}. 
Since $\cP_{3,3}^+ \cong \cP_{3,6}^e$ by the correspondence 
$f(x,y,z) \to f(x^2,y^2,z^2)$, 
our results characterize $\cE(\cP_{3,6}^e)$, 
and prove the following (see \Ttanc): 

\proclaim{Corollary 1.8} 
{\sl $\cE(\cP_{3,6}^e) \subset \cE(\cP_{3,6})$. }
\endproclaim

\medskip

Here, we sketch our idea of proof of \Ttae. 
To classify $f \in \cE(\cP_{3,3}^+)$, 
we observe the complex cubic curve $C_{\C} := V_{\C}(f) \subset \P_{\C}^2$ 
and the real cubic curve $C_{\R} := V_{\R}(f) \subset \P_{\R}^2$. 

If $C_{\C}$ is reducible, classification is easy (\Ttaaa---\Ttadan). 

Consider the case $C_{\C}$ is irreducible. 
Then $C_{\C}$ is a rational curve with a singular point $P_0$ (\Ttaea). 
If $P_0$ is a node, it lies inside of $\P_+^2$ (\Ttafb---\Ttahn). 
After a suitable projective transformation, 
we may assume $P_0 = (1 \colon 1 \colon 1)$. 
If $P_0$ is a cusp, it lies on an edge of a triangle $\partial \P_+^2$ 
and it is not a vertex (\Ttaic, \Ttafbn). 
In this case, we may assume $P_0 = (0 \colon 1 \colon 1)$. 
In the both cases, $C_{\R} - \{P_0\}$ contacts to $\partial \P_+^2$ at 
some points $P_1$,$\ldots$, $P_r$.

We shall study a general theory of infinitely near zeros in \S2. 
As a monic polynomial in one variable is determined by its roots 
when we count the multiplicity exactly, 
an extremal ray $\R_+ \cdot f$ is essentially determined 
by $V_+(f)$ when we consider infinitely near zeros. 
By \Tqbae, $f \in \cE(\cP_{3,3}^+)$ satisfies 
\[\R_+ \cdot f = \cP_0 \cap \cP_1 \cap \cdots \cap \cP_r\]
where $\cP_i$ is a local cone or an infinitesimal local cone 
of $\cP_{3,3}^+$ at $P_i$ (see \Tqbc). 
$\cP_i$ has informationis on infinitely near points of $P_i$. 
This condition also gives a lower bound of $r$. 
On the other hand, at most two $P_i$ can exist 
on an edge of the triangle $\partial \P_+^2$. 
If two points $P_i$, $P_j$ exist on an edge, one of them must be a vertex. 
Thus $r \leq 4$. 
Our classification essentially depends on this result. 

In \S 3.3, we study the case that $P_0 = (0 \colon 1 \colon 1)$ is a cusp. 
In this case, we can easily prove $f = c \frh_{pq}$ (see \Ttcb). 

In \S 3.4, we study the case that $P_0 = (1 \colon 1 \colon 1)$ is a node. 
In this case, there are several types 
of configurations of $P_1$,$\ldots$, $P_r$. 
After studying each case using \Ttac \ and \Ttadn, 
we conclude that $f$ can be represented by $\frf_{pqr}$ or $\frg_{pq}$. 

\medskip

The idea of infinitely near zeros is also useful to reform \Ttaf. 
If $f \in \cE(\cP_{3,6}) - \Sigma_{3,6}$ is not always exposed, 
$V_{\R}(f)$ consists of just 10 points including infinitely near points. 
\Tqbae \ states that 
extremal elements of PSD cones $\cP_{n,2d}$ or $\cP_{n,d}^+$ are 
determined by their equality conditions including infinitely near points. 
In \cite{RefR}, properties of $f \in \cE(\cP_{3,6}) - \Sigma_{3,6}$ with 
$\# V_{\R}(f) \leq 9$ are studied 
using the idea of `divisor' instead of `infinitely near points'. 
When we treat $V_{\R}(f)$ as a divisor, 
a notion of multiplicity is included in it, 
and works fairly well for $(n$, $2d)=(3$, $6)$. 
But it will not determine $f \in \cE(\cP_{n,2d})$ with real zeros of 
high multiplicity. 
Notion of `infinitely near points' gives more precise information 
about degeneration. 

We also mention that the $\R$-scheme structure of $\P_{\R}^2$ is not unique. 
We can blow up $\P_{\R}^2$ at conjugate imaginal points. 
So, we must be careful to treat `divisors' on $\P_{\R}^2$. 
Note that there also exist infinitely many real projective surfaces $X$ such 
that $\P_+^2 \subset X$. 
We often need to get rid of information outside of $\P_+^2$, 
and there shoud be the unique structure sheaf $\cR_{\P_+^2}$. 
So, in this article, we use the notion of semialgebraic variety. 

\section{Infinitely near zeros}

One of the key idea to prove \Ttae \ is to introduce the notion 
of infinitely near zeros. 
Max Noether introduced the notion of infinitely near points 
on algebraic surfaces. 
The notion of infinitely near zeros of inequalities is similar one. 
We need blowing ups to treat infinitely near zeros. 
So, we should generalize some notions. 
The notion of semialgebraic varieties is introduced in \cite{RefAa}, and 
the precise properties of semialgebraic variety are explained 
in \cite[\S 5]{RefAd}. 
But we don't need deep understanding for semialgebraic varieties in this 
article. 
We present here minimum definitions. 

\proclaim{Definition 2.1}{\rm (Semialgebraic variety)} 
A locally ringed space $(A$, $\cR_A)$ is called 
{\it semialgebraic variety}, if there exists 
a real algebraic variety $(X$, $\cR_X)$ in the sence of \cite{RefBCR} 
which satisfies the following: \par
{\parindent=20pt
\Item{(1)} The ring $\cR_X(X)$ is an integral domain such that 
$\Krulldim \cR_X(X) = \dim X$. 
We denote $\Rat(X) := Q(\cR_X(X))$ (the field of fractions). 
\Item{(2)} There exists an injective morphism $\iota \colon (A$, $\cR_A) 
 \lto (X$, $\cR_X)$ as locally ringed spaces, 
and $\iota$ induces a homeomorphism $A \to \iota(A)$ 
with respect to the Euclidean topology. 
Moreover, $\iota(A)$ is a semialgebraic subset of $X$ 
such that $\Zar_X(\iota(A)) = X$, i.e. 
the Zariski closure of $\iota(A)$ in $X$ agrees with $X$. 
\Item{(3)} The induced map
 $\iota_P^* \colon \cR_{X,\iota(P)} \lto \cR_{A,P}$ 
is an isomorphism for every $P \in A$, and 
$$\cR_A(U) = \bigcap_{P \in U} \iota_P^*\big(\cR_{X,\iota(P)}\big) 
  \subset \iota^* \Rat(X)$$
for any non-empty Euclidean open subset $U \subset A$. 

}
\endproclaim

The above definition is based on the fact that $\cR_X$ is generated by its 
global section (see Corollary 5.13 of \cite{RefAd}). 
On a semialgebraic variety $A$, we use the Euclidean topology and 
the Zariski topology induced from $X$. 
Many terminologies for semialgebraic varieties can be defined 
by the similar way as complex algebraic varieties. 

\def\Tqbv{Definition 2.2} 
\proclaim{Definition 2.2}{\rm (Signed linear system)} 
{\rm Let $(A$, $\cR_A)$ be a semialgebraic variety, 
and $\cC^0_A$ be the sheaf of germs of real continuous functions on A. \par
{\parindent=20pt
\Item{(1)} Let $\cI$ be an invertible $\cR_A$-sheaf. 
$\cI$ is called a {\it signed invertible sheaf} on $A$ if 
{\parindent=40pt
\Item{(i)} there exists $\cC^0_A$-invertible sheaf $\cJ$ 
such that $\cI \otimes_{\cR_A} \cC^0_A = \cJ \otimes_{\cC^0_A} \cJ$, and 
\Item{(ii)} there exists $e \in \cJ(A)$ such 
that $e^2 \in \cI(A)$ and $\cI(A) = \cR_A(A) \cdot e^2$. 

}
Then, for $f \in H^0(A$, $\cI)$, 
there exists $g \in H^0(A$, $\cR_A)$ such that $f = g e^2$. 
We define $\sign(f(P)) \in \{0$, $\pm 1\}$ by $\sign(f(P)) = \sign(g(P))$ 
for $P \in A$. 
\Item{(2)} Let $\cI$ be a signed invertible $\cR_A$-sheaf. 
A finite dimensional vector subspace $\cH \subset H^0(A$, $\cI)$ is 
called a {\it signed linear system} on $A$. 
For $f \in \cH$, we say $f$ is {\it PSD} on $A$ if 
$f(P) \geq 0$ for all $P \in A$. 

}
\endproclaim

\proclaim{Example 2.3} 
Let $A = \P_+^{n-1} \subset \P_{\R}^{n-1} = X$ and $e := \sqrt{x_1}^d$. 
Then $\cR_A(d) = \cR_A \cdot e^2$. 
Thus $\cR_A(d)$ is a signed invertible $\cR_A$-sheaf, 
and $\cH_{n,d} \subset H^0(A$, $\cR_A(d))$ is a signed linear system 
on $\P_+^n$, 

Similarly, $\cH_{n,2d}$ is a signed linear system on $\P_{\R}^{n-1}$. 
\endproclaim

Old \cite[Definition 1.7]{RefAa} sould be replaced by 
the above Definition 2.2. 
Otherwise, the following proposition does not hold. 

\proclaim{Proposition 2.4}
{\sl Let $(A$, $\cR_A)$ and $(B$, $\cR_B)$ be semialgebraic varieties 
and $\varphi \colon B \to A$ be a morphism. 
If $\cH$ is is a signed linear system on $A$, 
then $\varphi^*\cH := \big\{ \varphi^*(f) := f \circ \varphi$ 
$\big|$ $f \in \cH \big\}$ is a signed linear system on $B$. }
\endproclaim

\Proof
Let $\cI$, $\cJ$ and $e$ be same as in \Tqbv. 
Then $\cI_B := \varphi^* \cI \otimes_{\varphi^*\cR_A} \cR_B$, 
$\cJ_B := \varphi^* \cJ \otimes_{\varphi^*\cC^0_A} \cC^0_B$ 
and $e_B := \varphi^*(e)$ satisfy conditions 
so that $\varphi^*\cH$ is a signed linear system on $B$. 
\end{proof}

\proclaim{Definition 2.5}{\rm (PSD cone)} 
Let $(A$, $\cR_A)$ be a semialgebraic variety, 
and $\cH$ be a signed linear system on $A$. The cone 
\[\cP = \cP(A, \, \cH) := \big\{f \in \cH \; \big| \; 
      \hbox{$f$ is PSD on $A$}\big\}\]
is called the {\it PSD cone} on $A$ in $\cH$. 
\endproclaim

We can represent $\cP_{n,d}^+ = \cP(\P_+^{n-1}$, $\cH_{n,d})$, 
and $\cP_{n,2d} = \cP(\P_{\R}^{n-1}$, $\cH_{n,2d})$. 
In this article, we only treat the cases that $\cP(A$, $\cH)$ contains 
no lines. 
So, when we say '$\cP$ is a PSD cone', 
we always assume that $\cP$ contains no lines. 

\def\Tqba{Definition 2.6} 
\proclaim{Definition 2.6} 
{\rm Let $A$ be a non-singular semialgebraic variety, 
$\cH$ be a signed linear system on $A$, 
and $\cP = \cP(A$, $\cH)$ be a PSD cone. \par
{\parindent=20pt
\Item{(1)} Take $P \in A$ and $f \in \cH$. 
Assume that $f$ can be represented as $f = g e^2$ as \Tqbv. 
If we take a suitable open subset $U \subset A$ 
and an analytic coordinate system $(x_1$,$\ldots$, $x_n)$ 
on $U$ whose origin is $P$, we can regard 
$g|_U \in \hat{\cR}_{A,P} := \R[[x_1$,$\ldots$, $x_n]]$. 
Let $\mathfrak{m}$ be the maximal ideal of $\hat{\cR}_{A,P}$ corresponding 
to the point $P$. 
The {\it multiplicity} of $f$ at $P$ is defined as
\[\mult_P f 
  := \sup \big\{ d \geq 0 \; \big| \; g|_U \in \mathfrak{m}^d \big\}.\]
\Item{(2)} For $x \in A$, we put 
\[\mult_x \cP := \min_{h \in \cP-\{0\}} \mult_x h.\]
Assume that $\dim \cP \geq 2$, 
where $\dim \cP$ implies the dimension of $\cP$ as a semialgebraic 
variety which agree with the dimension as a convex cone or as a manifold. 
For $0 \ne f \in \cP$, we put 
\[\cZ_f(\cP) := \big\{ P \in A \; \big| \; 
  \mult_P f > \mult_P \cP \big\}.\]
\Item{(3)} Take a point $P \in A$ and put 
\[\cP_P := \big\{ g \in \cP \; \big| \; 
  \mult_P g > \mult_P \cP\big\}.\]
If $\cP_P \ne \{0\}$, then $\cP_P$ is called 
the {\it local cone} of $\cP$ at $P$. 

}
\endproclaim

When we consider a vector subspace $\cH_P := \big\{ g \in \cH$ $\big|$ 
$\mult_P g > \mult_P \cP\big\}$, then $\cP_P = \cP(A$, $\cH_P)$. 
Thus $\cP_P$ is a semialgebraic closed convex cone. 

In \cite{RefAa}, the author gave a different definition of $\cP_P$ as 
\[\cP_P := \big\{ f \in \cP \; \big| \; 
   \hbox{$f(P) = 0$}\big\}.\]
If $P \notin \Bs \cH$, this agrees with $\cP_P$ in the above definition (3). 
In the case $P \in \Bs \cH$, it seems that the above Definiton 2.6 
of $\cP_P$ works well. 

The word `local cone' is used in \cite{RefSchu} to 
describe the cone of locally non-negative forms at $P$. 
Please don't confuse with the above definition. 

\def\Tqbw{Proposition 2.7} 
\proclaim{Proposition 2.7} 
{\sl Let $A$, $\cP$ be same as 
in \Tqba \ where $\dim \cP \geq 2$. }\par
{\parindent=20pt
\Item{(1)} {\sl There exists $h_0 \in \cP$ such that 
$\mult_x h_0 = \mult_x \cP$ for all $x \in A$. }
\Item{(2)} {\sl Assume that $0 \ne f \in \cP_a$ and $f = g + h$, 
where $g$, $h \in \cP$ and $a \in A$. 
Then $g$, $h \in \cP_a$. }
\endproclaim

\Proof
(1) For $a \in A$, there exists $h_a \in \cP$ such that 
$\mult_a h_a = \mult_a \cP$. 
Note that $\mult_x g$ and 
$\mult_x \cP$ are upper semicontinuous functions on $x \in A$ with 
respect to Zariski topology (see \cite[II Exercise 5.8]{RefHart}). 
For any point $a \in A$, there exists an 
open neighborhood $a \in U \subset A$ and $h \in \cP$ 
such that $\mult_x h = \mult_x \cP$ for all $x \in U$. 
Since $\mult_x h_a \geq \mult_x \cP$, $h_a/h$ is holomorphic on $U$ 
and is upper semicontinuous. 
Since $\mult_x h_a - \mult_x \cP = \mult_x h_a/h \in \Z$, 
this is is also a Zariski upper semicontinuous functions on $x \in A$, 
and whose minimum value is equal to $0$. 
Thus 
$U_a := \big\{ x \in A$ $\big|$ $\mult_x h_a = \mult_x \cP \big\}$. 
is a Zariski open subset of $A$. 
Since $A$ is quasi-compact with respect to Zariski topology, 
we can choose $a_1$,$\ldots$, $a_r \in A$ such 
that $U_{a_1} \cup \cdots \cup U_{a_r} = A$. 
Let $h_0 = h_{a_1} + \cdots + h_{a_r}$. Then 
\[\mult_x h_0 = \inf_{1 \leq i \leq r} \mult_x h_{a_i} = \mult_x \cP\]
for all $x \in A$. 

\smallskip

(2) Take a general $h_0 \in \Int(\cP)$ 
such that $\mult_x h_0 = \mult_x \cP$ for all $x \in A$, 
where the word `general' is used 
in the sense \cite[\S 7.9 {\bf a}]{RefIItaka}. 
Then $f_1 := f/h_0$, $g_1 := g/h_0$ and $h_1 := h/h_0$ 
are holomorphic functions on $A$ (i.e. $f_1$, $g_1$, $h_1 \in 
H^0(A$, $\cR_A)$), for $f_1$, $g_1$, $h_1$ have no poles on $A$ by (1). 
Moreover, $f_1$, $g_1$ and $h_1$ are PSD on $A$. 
Since $0 = f_1(a) = g_1(a) + h_1(a)$, $g_1(a) \geq 0$ and $h_1(a) \geq 0$, 
we have $g_1(a) = h_1(a) = 0$. 
Thus $\mult_a g > \mult_a \cP$ and $\mult_a h > \mult_a \cP$. 
This implies $g$, $h \in \cP_a$. 
\end{proof}

\def\Tqbb{Proposition 2.8} 
\proclaim{Proposition 2.8} 
{\sl Let $A$, $\cP$ be same as in \Tqba \ with $\dim \cP \geq 2$. 
Moreover, we assume that $A$ is compact 
with respect to the Euclidean topology. 
Take $f \in \cE(\cP)$. Then, } \par
{\parindent=20pt
\Item{(1)} {\sl $\cZ_f(\cP) \ne \emptyset$. } 
\Item{(2)} {\sl If $a \in \cZ_f(\cP)$, 
then $\cE(\cP_a) = \cP_a \cap \cE(\cP)$. }

}
\endproclaim

\Proof
(1) Assume that $\cZ_f(\cP) = \emptyset$. 
Then, $f \in \cP$ satisfies 
$\mult_x f = \mult_x \cP$ for all $x \in A$. 
Moreover, a general $h_0 \in \Int(\cP) - \R_+ \cdot f$ satisfy 
$\mult_x h_0 = \mult_x \cP$ for all $x \in A$. 
We can regard $g := f/h_0$ as a holomorphic function on $A$. 
Since $A$ is compact, 
$\displaystyle \varepsilon := \inf_{x \in A} g(x) > 0$. 
Then $h := f - \varepsilon h_0 \in \cP - \R_+ \cdot f$. 
Thus $f = h + \varepsilon h_0$ is not extremal in $\cP$. 

\smallskip

(2) $\cE(\cP_a) \supset \cP_a \cap \cE(\cP)$ is 
trivial. We shall show $\cE(\cP_a) \subset \cE(\cP)$. 
Take $g \in \cE(\cP_a)$. 
Assume that $g = h_1 + h_2$ for $h_1$, $h_2 \in \cP - \R_+ \cdot g$. 
Then $h_1$, $h_2 \in \cP_a$ by \Tqbw. 
So, $g \notin \cE(\cP_a)$. A contradiction. 
Thus $g \in \cE(\cP)$, 
and $\cE(\cP_a) \subset \cE(\cP)$. 
\end{proof}

Note that $\cE(\cQ) = \cQ \cap \cE(\cP)$ does not hold 
for PSD cones $\cQ \subset \cP$, in general. 
For example, $\cP_{3,4} \subset \cP_{3,4}^+$ but 
$\cE(\cP_{3,4}) \not\subset \cE(\cP_{3,4}^+)$. 

The following definition is an analogue of 
a resolution of the base locus of a linear system by 
a sequence of blowing ups (see \cite{RefHironaka}). 
Words and symbols are based on algebraic geometry. 

\def\Tqbc{Definition 2.9} 
\proclaim{Definition 2.9}{\rm (infinitesimal local cone)} 
{\rm Let $A$ be a non-singular semialgebraic variety 
which is compact with respect to the Euclidean topology. 
Let $\cP = \cP(A$, $\cH)$ be a PSD cone 
with $\dim \cP \geq 2$. 
Fix $f \in \cE(\cP)$. 
Then $\cZ_f(\cP) \ne \emptyset$. 

\smallskip

(1) Take $a \in \cZ_f(\cP)$. 
Assume that $\dim \cP_a \geq 2$. 
Put $A_0 := A$, $a_0 := a$, $f_0 := f$ and $\cL_0 := \cP_a$. 
Then $f_0 \in \cE(\cL_0)$ by \Tqbb. 

Inductively, we shall define $A_i$, $a_i$, $\cL_i$ for $i \geq 0$, 
and $\psi_i \colon A_i \to A_{i-1}$ for $i \geq 1$. 
Now fix $i \in \N$, and assume that $A_j$, $a_j$ and $\cL_j$ 
are defined for all $0 \leq j < i$. 
Put $\overline{\psi}_j := \psi_1 \circ \cdots \circ \psi_j : A_j \lto A$ 
whenever $\psi_1$,$\ldots$, $\psi_j$ will be defined. 
In this process, we assume that $\psi_j(a_j) = a_{j-1}$, 
$f_j := \overline{\psi}_j^* (f) = f \circ \overline{\psi}_j \in \cE(\cL_j)$ 
and $a_j \in \cZ_{f_j}(\cL_j)$ for $0 \leq j < i$. 
Consider 
\[\cL_i' := \big\{ g \in \cL_{i-1} \; \big| \; 
  \mult_{a_{i-1}} g > \mult_{a_{i-1}} \cL_{i-1}\big\}.\]
We divide into two cases. 

\smallskip

(1-i) The case $0 \ne \cL_i' \subsetne \cL_{i-1}$ 
and $f_{i-1} \in \cL_i'$. 

Then, let $A_i := A_{i-1}$, $\psi_i \colon A_i \lto A_{i-1}$ be 
the identity map, $a_i := a_{i-1}$ and 
we put $\cL_i := \cL_i'$. 
Note that $a_i \in \cZ_{f_i}(\cL_i)$. 
Now we repeat the process increasing $i$. 

\smallskip

(1-ii) The case $\cL_i' = 0$ or $\cL_i' = \cL_{i-1}$ or 
$f_{i-1} \notin \cL_i'$. 

Then, let $\psi_i \colon A_i \lto A_{i-1}$ be the blowing up 
of $A_{i-1}$ at the point $a_{i-1} \in A_{i-1}$. 
$\cH_i := \overline{\psi}_i^* \cH$ is 
a signed linear system on $A_i$ and 
$\overline{\psi}_i^* \cP 
 := \big\{ g \circ \overline{\psi}_i$ $\big|$ $g \in \cP \big\}$ 
satisfies $\overline{\psi}_i^* \cP 
  = \cP(A_i$, $\cH_i)$.  
Note that $\dim \psi_i^* \cL_{i-1} = \dim \cL_{i-1} \geq 2$. 

\smallskip

(1-ii-a) Consider the case that we can find 
$a_i \in \cZ_{f_i}(\psi_i^* \cL_{i-1})$ such that 
$\psi_i(a_i) = a_{i-1}$ where $f_i := \overline{\psi}_i^* f$. 
Let $\cL_i$ be the local cone of $\psi_i^* \cL_{i-1}$ at the 
point $a_i$. 
Note that $f_i \in \cE(\cL_i)$, 
because $\cE(\cL_i) 
 = \cL_i \cap \cE(\psi_i^* \cL_{i-1})$. 
Then, we repeat the process increasing $i$. 

\smallskip

(1-ii-b) Termination of the process. 

Since $\dim \cP > \dim \cL_0 > \dim \cL_1 > \cdots$, 
there exists $l \in \N$ such that 
$\big\{a \in \cZ_{f_{l+1}}(\psi_{l+1}^* \cL_l)$ $\big|$ 
$\psi_{l+1}(a) = a_l \big\} = \emptyset$. 
Then, we stop to repeat the process. 
We say $a_1$, $a_2$,$\ldots$, $a_l$ is a sequence of zeros of $f$ 
{\it infinitely near} to $a$. 
Each $a_i$ ($1 \leq i \leq l$) is called a zero of $f$ infinitely near to $a$ 
in $\cP$. 
The convex cone $\overline{\psi_l}(\cL_l)$ is called an 
{\it infinitesimal local cone} of $\cP$ at $a_l$ or at $a$ 
with respect to $f$. 

Assume that $f(a)=0$ and 
$f$ has only finitely many zeros $b_1$,$\ldots$, $b_N$ 
infinitely near to $a$ in $\cP$. 
Then, we define $\length_a f := N+1$. 
If $\dim \cP_a = 1$ or there exists no 
$a_1 \in \cZ_{\psi_1^* f}(\psi_1^* \cP_a)$ such that 
$\psi_1(a_1) = a$, then we put $\length_a f := 1$. 
If $f(a) \ne 0$, then we put $\length_a f := 0$. }
\endproclaim

\def\Tqbf{Lemma 2.10} 
\proclaim{Lemma 2.10} 
{\sl Let $A$, $\cP$ and $f \in \cE(\cP)$ be as in \Tqbc. 
Assume that $a$, $b \in \cZ_f(\cP)$ with $a \ne b$. 
Let $\cQ := \cP_b$, and assume that $\dim \cQ \geq 2$. 
If $\cL$ is a local cone or an infinitesimal local cone of $\cQ$ at 
a point $a \in A$, then there exists a local cone or 
an infinitesimal local cone $\overline\cL$ of $\cP$ at $a$ 
such that $\overline\cL \cap \cQ = \cL$. }
\endproclaim

\Proof
If $\cL$ is an infinitesimal local cone of $\cQ$, 
we take a sequences $\{A_i\}_{i=0}^l$, $\{a_i\}_{i=0}^l$, 
$\{\psi_i\}_{i=1}^l$ and $\{\cL_i\}_{i=0}^l$ as in \Tqbc, 
such that $\cL_0 = \cQ_a$, 
$\overline{\psi_l}(\cL_l) = \cL$, $A_0 = A$ and $f_0 = f$. 

In the case that $\cL = \cQ_a$ is a local cone of $\cQ$, 
we put $l=0$ at the above sequences. 

Put $\overline\cL_0 := \cP_a$. 
We need to find convex cones $\overline\cL_i$ on $A_i$ such that 
$\overline\cL_i \cap \overline{\psi_i}^* \cQ = \cL_i$ 
and that $\overline\cL_i$ is a local cone or an infinitesimal 
local cone of $\psi_i^* \cL_{i-1}$. 
To construct $\overline\cL_i$, we may need to refine the sequences. 
In the process of the refinement, 
we always put $A_i^0 = A_i^1 = \cdots = A_i^{k_i} = A_i$ ($0 \leq i \leq l$), 
and $\psi_i^j = \id \colon A_i^{j+1} \to A_i^j$ 
($0 \leq j < k_i$, $0 \leq i \leq l$) 
and $\psi_{i+1}^{k_i} = \psi_{i+1} \colon A_{i+1}^0 \lto A_i^{k_i}$ 
($0 \leq i < l$). 

\smallskip

(1) If $1 \leq i \leq l$ and $\overline\cL_{i-1}$ is already defined 
so that $\overline{\cL_{i-1}} \cap \overline{\psi_{i-1}}^* \cQ = \cL_{i-1}$, 
we put $\cL_i^0 := \big(\psi_{i-1}^* \overline\cL_{i-1}\big)_{a_i}$. 

It is easy to see that 
$\cL_i^0 \cap \overline{\psi_i}^*\cQ \supset \cL_i$. 

(2) If $\cL_i^0 \cap \overline{\psi_i}^*\cQ = \cL_i$, 
  then we put $\overline\cL_i := \cL_i^0$, $k_i := 0$, 
and we don't need a refinement. 

(3) Consider the case 
$\cL_i^0 \cap \overline{\psi_i}^*\cQ \supsetne \cL_i$. 

This can happen only if $\mult_{a_i} \cL_i^0 < \mult_{a_i} \cL_i$. 
Let $A_i^1 := A_i$ and $\psi_i^1 = \id \colon A_i^1 \to A_i^0$. 
Put $\cL_i^1 := (\cL_i^0)_{a_i}$. 
Since $\mult_{a_i} \cL_i \leq \mult_{a_i} f_i$, 
we have $f_i \in \cL_i^1 \subsetne \cL_i^0$. 
Note that $\mult_{a_i} \cL_i^0 < \mult_{a_i} \cL_i^1 
  \leq \mult_{a_i} \cL_i \leq \mult_{a_i} f_i$. 
Thus, $\cL_i^1 \cap \overline{\psi_i}^*\cQ \supset \cL_i$. 

(4) If $\cL_i^1 \cap \overline{\psi_i}^*\cQ = \cL_i$, 
  then we put $\overline\cL_i := \cL_i^1$, $k_i := 1$, 
and we stop this refinement. 

(5) If $\cL_i^1 \cap \overline{\psi_i}^*\cQ \supsetne \cL_i$, 
Put $\cL_i^2 := (\cL_i^1)_{a_i}$. 
Repeat this process till 
$\cL_i^{k_i} \cap \overline{\psi_i}^*\cQ = \cL_i$. 
Then we put $\overline\cL_i := \cL_i^{k_i}$. 

When $i=l$, we put 
$\overline\cL := \overline{\psi_l}(\overline\cL_l)$. 
Then  $\overline\cL \cap \cQ = \cL$. 
\end{proof}

\def\Tqbae{Theorem 2.11} 
\proclaim{Theorem 2.11} 
{\sl Let $A$, $\cP$ and $f \in \cE(\cP)$ be as in \Tqbc. 
Assume that $\dim \cP \geq 2$. 
Then, there exists points $P_1$,$\ldots$, $P_r \in A$ 
(not always distinct), and local cones or infinitesimal local cones 
$\cP_1$,$\ldots$, $\cP_r \subset \cP$ with respect to $f$ 
which satisfy } \par
{\parindent=20pt
\Item{(1)} $\cP_1 \cap \cdots \cap \cP_r = \R_+ \cdot f$. 
\Item{(2)} {\sl $\cP_i$ is the local cone $\cP_{P_i}$ or 
an infinitesimal local cone of $\cP$ at $P_i \in A$ 
with respect to $f$ for $i=1$,$\ldots$, $r$. }

}
\endproclaim

\Proof
We prove by induction on $\dim \cP$. 
Take $c \in \cZ_f(\cP)$, and put $\cQ := \cP_c$. 

If $\dim \cP = 2$, then $\cP_c$ is a local cone and $\cP_c = \R_+ \cdot f$ 
(note that we assume $\cP$ contains no lines.) 

Consider the case $\dim \cP \geq 3$. 

Since $\dim \cQ < \dim \cP$, 
there exists points $P_1$,$\ldots$, $P_r \in A$, and 
local cones or infinitesimal local cones 
$\cQ_1$,$\ldots$, $\cQ_r$ of $\cQ$ which satisfy 
$\cQ_1 \cap \cdots \cap \cQ_r = \R_+ \cdot f$, 
and that $\cQ_j$ is the local cone $\cQ_{P_j}$ or 
an infinitesimal local cone of $\cQ$ at $P_j$. 
Then there exists a local cone or 
an infinitesimal local cone $\cP_j$ of $\cP$ at $P_j$ 
such that $\cP_j \cap \cQ = \cQ_j$ by the above lemma. 
Thus
\[\cP_c \cap \cP_1 \cap \cdots \cap \cP_r 
  = \cQ_1 \cap \cdots \cap \cQ_r = \R_+ \cdot f.\]
If there exists $1 \leq i \leq r$ such that $\cP_i \subset \cP_c$, 
we have $\cP_1 \cap \cdots \cap \cP_r = \R_+ \cdot f$. 
Otherwise, put $P_{r+1} := c$ and $\cP_{r+1} := \cP_c$. 
Then $\cP_1 \cap \cdots \cap \cP_{r+1} = \R_+ \cdot f$. 
\end{proof}

\def\Tqbg{Definition 2.12} 
\proclaim{Definition 2.12} 
{\rm Let $A$ be a non-singular semialgebraic variety 
which is compact with respect to the Euclidean topology. 
Let $\cP = \cP(A$, $\cH)$ be a PSD cone. 
Take $a \in A$ and put $\cL_0 := \cP_a$, $A_0 := A$, $a_0 = a$. 
Assume that $\dim \cL_0 \geq 2$. 

Let $\{A_i\}_{i=0}^l$, $\{a_i\}_{i=0}^l$, $\{\psi_i\}_{i=1}^l$, 
$\{\cL_i\}_{i=0}^l$ be the sequence such that : \par
{\parindent=20pt
\Item{(1)} $\psi_i \colon A_i \lto A_{i-1}$ is the blowing up 
of $A_{i-1}$ at the point $a_{i-1} \in A_{i-1}$ or 
the identity map ($A_i = A_{i-1}$, $\psi_i = \id$) ($1 \leq i \leq l$). 
\Item{(2)} $\psi_i(a_i) = a_{i-1}$ ($1 \leq i \leq l$). 
\Item{(3)} $\cL_i$ is the local cone of $\psi_i^* \cL_{i-1}$ at 
the point $a_i \in A_i$ ($1 \leq i \leq l$). 
\Item{(4)} $\dim \cL_i \geq 2$ for $1 \leq i < l$ and 
$\dim \cL_l \geq 1$. 

}
Put $\overline{\psi}_i := \psi_1 \circ \cdots \circ \psi_i : A_i \lto A$ 
($1 \leq i \leq l$). 
Then $\overline{\psi}_l(\cL_l)$ is called an 
{\it infinitesimal local cone} of $\cP$ at $a$. }
\endproclaim

\def\Tqbh{Proposition 2.13} 
\proclaim{Proposition 2.13} 
{\sl Let $A$ be a non-singular compact semialgebraic variety. 
$\cP$ be a PSD cone on $A$, and $f \in \cP$. 
Assume that $\dim \cP \geq 2$, 
and there exists local cones or infinitesimal local cones 
$\cP_1$,$\ldots$, $\cP_r \subset \cP$ such that 
$\cP_1 \cap \cdots \cap \cP_r = \R_+ \cdot f$. 
Then, $f \in \cE(\cP)$. }
\endproclaim

\Proof
Assume that there exists $g$, $h \in \cP - \R_+ \cdot f$ 
such that $f = g + h$. 
Fix $1 \leq k \leq r$. 

\smallskip

(1) We shall prove that $g$, $h \in \cP_k$. 

(1-i) If $\cP_k$ is a local cone, (1) follow from \Tqbw. 

(1-ii) Consider the case that $\cP_k$ is 
an infinitesimal local cone of $\cP$ at $a$. 
Take sequences $\{A_i\}_{i=0}^l$, $\{a_i\}_{i=0}^l$, $\{\psi_i\}_{i=1}^l$, 
$\{\cL_i\}_{i=0}^l$ 
with $\overline{\psi}_l(\cL_l) = \cP_k$ as in the above definition. 
Formally, put $\overline{\psi}_0 := \id \colon A_0 \to A_0$. 
Put $f_i := \overline{\psi}_i^* (f)$, $g_i := \overline{\psi}_i^* (g)$ and 
$h_i := \overline{\psi}_i^* (h)$ ($0 \leq i \leq l$). 
Then $f_i = g_i + h_i \in \cL_i$. 
We shall show that $g_i$, $h_i \in \cL_i$ by induction on $i$. 

If $i=0$, we have $g_0$, $h_0 \in \cP_a = \cL_0$ by \Tqbw. 

Assume that $i \geq 1$ and $g_{i-1}$, $h_{i-1} \in \cL_{i-1}$. 
$\cL_i$ is the local cone of $\psi_i^* \cL_{i-1}$ at 
the point $a_i \in A_i$. 
Thus, we have $g_i$, $h_i \in \cL_i$, by \Tqbw. 

\smallskip

(2) By (1), we have $g$, $h \in \cP_1 \cap \cdots \cap \cP_r
   = \R_+ \cdot f$. 
Thus, $f \in \cE(\cP)$. 
\end{proof}

Let $f \in \cH_{3,d}$, $g \in \cH_{3,e}$ and $P \in \P_{\C}^2$. 
Take a local coordinate system $(x$, $y)$ on an affine open subset $U 
\subset \P_{\C}^2$ whose origin is $P$. 
We consider $f$, $g \in \C[[x$, $y]]$ and 
we denote the local intersection number of $V_{\C}(f)$ and $V_{\C}(g)$ at 
$P$ by $I_P(f,g) := \dim_{\C} \C[[x,y]]/(f,g)$. 
The intersection number of $V_{\C}(f)$ and $V_{\C}(g)$ is denoted by $I(f,g)$. 

\def\Tqbai{Example 2.14} 
\proclaim{Example 2.14} 
{\rm Consider the case $A = \P_+^2$ and 
$\cP = \cP_{3,d}^+$ ($d \geq 2$). 
We denote the homogeneous coordinate system of $A$ by 
$(x_0 \colon x_1 \colon x_2)$. 

\smallskip

(1) Consider $\cP_P$ when $P = (1 \colon p \colon q) \in \P_+^2$. 
Put $x := (x_1 - p x_0)/x_0$ and $y := (x_2 - q x_0)/x_0$. 
Take an arbitrary $f \in \cP_P$. Then 
$f(x,y) = ax^2 + 2bxy + cy^2 + (\hbox{\rm higher terms})$ such that 
$\displaystyle C_f := 
\left(\begin{matrix}a & b \\ b & c \end{matrix}\right)$ is 
positive semidefinite. 
Let $\psi_1 \colon A_1 \to A$ be the blowing up at $P$ and put $t := x/y$.
Then $f_1 := \psi_1^* (f) = y^2(at^2+2bt+c +\cdots)$ and 
$\mult_Q f_1 \geq 2$ for all $Q \in \psi_1^{-1}(P)$. In this case 
\begin{align*}
 \cP_P
 & = \big\{ f(x,y) \in \cP \; \big| \; 
      \hbox{$f(P)=f_x(P)=f_y(P)=0$} \big\} \\
 & = \big\{ F(x_0,x_1,x_2) \in \cP \; \big| \; 
      \hbox{$F(P)=F_{x_1}(P)=F_{x_2}(P)=0$} \big\} 
\end{align*}
where $f_x = \partial f(x,y)/\partial x$, 
$F_{x_0} = \partial F(x_0,x_1,x_2)/\partial {x_0}$ and so on. 

If $C_f$ is a positive definite matrix, 
Then $\mult_Q f_1 = 2$ for all $Q \in \psi_1^{-1}(P)$, 
and $\cZ_f(\psi_1^* \cP_P) = \emptyset$. 
Then $\length_P f = 1$. 

\smallskip

(2) Assume that $d \geq 4$ and 
the leading term of $f(x$, $y)$ is equal to $x^2+y^4$. 
Then $f_1 = y^2(t^2+y^2+\cdots)$ and 
$\cZ_f(\psi_1^* \cP_P)$ consists of a single point $P_1$ 
defined by $(t$, $y)=(0$, $0)$. 
Let $\psi_2 \colon A_2 \to A_1$ be the blowing up 
at $P_1$ and put $t_2 := x/t$.
Then $f_2 := \psi_2^* f_1 = y^4(t_2^2+1+\cdots)$. 
Thus $\length_P f = 2$. 
Let $\cP_{P_1}$ be the infinitesimal local cone of $\cP$ 
at $P_1$. Then 
\[\cP_{P_1} = \big\{ g(x,y) \in \cP \; \big| \; 
  \hbox{$g(P)=g_x(P)=g_y(P)=g_{yy}(P)=g_{xy}(P)=g_{yyy}(P)=0$} \big\}.\]

\smallskip

(3) Let $d := 3$, $P = (1 \colon p \colon 0) \in \P_+^2$ ($p>0$), 
$x := (x_1 - p x_0)/x_0$ and $y := x_2/x_0$. 
Note that $y \geq 0$ on $A$. 
Take a general $f \in \cP_P$. Then the leading term of $f$ is equal to 
$ay + bx^2$ with $a > 0$, $b>0$. 
Let $\psi_1 \colon A_1 \to A$ be the blowing up at $P$ and put $t := y/x$.
Then $f_1 = \psi_1^* (f) = x(at+bx + \cdots)$.
In this case 
\begin{align*}
 \cP_P 
 & = \big\{ f(x,y) \in \cP \; \big| \; \hbox{$f(P)=f_x(P)=0$} \big\} \\
 & = \big\{ F(x_0,x_1,x_2) \in \cP \; \big| \; 
      \hbox{$F(P)=F_{x_1}(P)=0$} \big\}.
\end{align*}
Fix $f \in \cP_P$. 
Then $\cZ_f(\psi_1^* \cP_P) = \emptyset$. 
Thus $\length_P f = 1$. 

\smallskip

(4) Let $P = (1 \colon 0 \colon 0) \in \P_+^2$, 
$x := x_1/x_0$ and $y := x_2/x_0$. 
Note that $x \geq 0$ and $y \geq 0$ on $A$. 
Assume that $f \in \cP_P$, $I_P(f$, $x) = 1$ and $I_P(f$, $y) = m$. Then 
$f(x$, $y) = ax + by^m + (\hbox{\rm higher terms})$ with $a>0$, $b>0$. 
Let $\psi_1 \colon A_1 \to A$ be the blowing up at $P$, 
and put $t_1 := x/y$. 
Then $f_1 = y(at_1 + by^{m-1} + \cdots)$. 
Let $P_1 \in A_1$ be the point defied by $(t_1$, $y) = (0$, $0)$. 

Similarly, let $\psi_i \colon A_i \to A_{i-1}$ be the blowing up at $P_i$, 
and put $t_i := t_{i-1}/y$. 
Let $P_i \in A_i$ be the point defied by $(t_i$, $y) = (0$, $0)$. 
Then $f_i = y^i(at_i + by^{m-i} + \cdots)$ ($i \leq m$). 
It is easy to see that $\length_P f = m$. 
If $m = 1$, then 
\begin{align*}
 \cP_P 
 & = \big\{ g(x,y) \in \cP \; \big| \; \hbox{$g(P)=0$} \big\} \\
 & = \big\{ F(x_0,x_1,x_2) \in \cP \; \big| \; 
      \hbox{$F(P)=0$} \big\}. 
\end{align*}
If $m = 2$, then 
\begin{align*}
 \cP_{P_1} 
 & = \big\{ g(x,y) \in \cP \; \big| \; 
        \hbox{$g(P)=g_y(P)=0$} \big\} \\
 & = \big\{ F(x_0,x_1,x_2) \in \cP \; \big| \; 
      \hbox{$F(P)=F_{x_2}(P)=0$} \big\}. 
\end{align*}
If $m = 3$, then 
\begin{align*}
 \cP_{P_2} 
 & = \big\{ g(x,y) \in \cP \; \big| \; 
        \hbox{$g(P)=g_y(P)=g_{yy}(P)=0$} \big\} \\
 & = \big\{ F(x_0,x_1,x_2) \in \cP \; \big| \; 
      \hbox{$F(P)=F_{x_2}(P)=F_{x_2x_2}(P)=0$} \big\}, 
\end{align*}
where $\cP_{P_i}$ be the infinitesimal local cone of $\cP$ 
at $P_i$. 
Thus $\length_P f = m$. 

\smallskip

(5) Let $P = (1 \colon p \colon 0) \in \P_+^2$ ($p>0$), 
$x := (x_1 - p x_0)/x_0$ and $y := x_2/x_0$. 
Assume that $f \in \cE(\cP_{3,3}^+)$ has the leading term 
$((y-1) + a x )^2 + b (x + c (y - 1))^3$ ($b \ne 0$). 
Let's determine all zeros of $f$ infinitely near to $P$. 
Put $v := (y-1) + a x$ and $u := x + c (y-1)$. 
Let $\psi_1 \colon A_1 \to A$ be the blowing up at $P$ and put $t := v/u$.
Then $f_1 = u^2((t^2+u) + \cdots)$. 
Let $P_1 \in A_1$ be the point defined by $(u$, $t)=(0$, $0)$. 
Since $\mult_{P_1} f_1 = 3 
 > 2 = \mult_{P_1} (\psi_1^* \cP_{3,3}^+)_{P_1}$, 
we need to blow up $\psi_2 \colon A_2 \to A_1$ at $P_1$. 
Put $s := u/t$.
Then $f_2 = u^2t((t+s) + \cdots)$. 
Let $P_2 \in A_2$ be the point defined by $(u$, $s)=(0$, $0)$. 
Then $\mult_{P_2} f_2 = 4 
 > 3 = \mult_{P_2} (\overline{\psi_2}^* \cP_{3,3}^+)_{P_2}$. 
There exists no more zero of $f$ infinitely near to $P$. 
In this case, 
\begin{align*}
 \cP_P 
  & = \big\{ F(x_0,x_1,x_2) \in \cP \; \big| \; 
      \hbox{$F(P)=F_{x_1}(P)=0$} \big\}, \\
 \cP_{P_1} 
  & = \big\{ g(u,v) \in \cP \; \big| \; 
        \hbox{$g(P)=g_u(P)=g_v(P)=0$} \big\} \\
  & = \big\{ F(x_0,x_1,x_2) \in \cP \; \big| \; 
      \hbox{$F(P)=F_{x_1}(P)=F_{x_2}(P)=0$.} \big\}, \\
 \cP_{P_2} 
  & = \big\{ g(u,v) \in \cP \; \big| \; 
        \hbox{$g(P)=g_u(P)=g_v(P)=g_{uu}(P)=g_{uv}(P)=0$} \big\} \\
  & = \left\{ F(x_0,x_1,x_2) \in \cP \ \left| \ 
    \vcenter{
      \hbox{$F(P)=F_{x_1}(P)=F_{x_2}(P)=0$,}
      \hbox{$F_{x_1x_1}(P)=F_{x_1x_2}(P)=0$}}\right.\right\}. 
\end{align*}
Thus $\length_P f = 3$. }
\endproclaim

\def\Tqbaib{Remark 2.15} 
\proclaim{Remark 2.15} 
{\rm (1) Let 
$F(x_0$, $x_1$, $x_2) \in \cH_{3,m}$ and $P \notin V_{\C}(x_0)$. 
Put $x := x_1/x_0$, $y := x_2/x_0$ and $f(x$, $y) := F(1$, $x$, $y)$. 
Since $f_x(x$, $y) = F_{x_1}(1$, $x$, $y)$, we have
\[\hbox{$f_x(P) = 0$ $\Longleftrightarrow$ $F_{x_1}(P) = 0$.}\]
Assume that $F(P)=0$ and $P = (1 \colon p \colon 0)$ with $p \ne 0$. 
Since $x_0 F_{x_0} + x_1 F_{x_1} + x_2 F_{x_2} = m F$, we have 
\[\hbox{$f_x(P) = 0$ $\Longleftrightarrow$ $F_{x_1}(P) = 0$
          $\Longleftrightarrow$ $F_{x_0}(P) = 0$,}\]
Assume that $F(P)=F_{x_1}(P)=0$ and $P = (1 \colon 0 \colon 0)$. 
Since $x_0 F_{x_0x_1} + x_1 F_{x_1x_1} + x_2 F_{x_1x_2} = (m-1) F_{x_1}$, 
$x_0 F_{x_0x_0} + x_1 F_{x_0x_1} + x_2 F_{x_0x_2} = (m-1) F_{x_0}$, we have 
\[\hbox{$f_{xx}(P) = 0$ $\Longleftrightarrow$ $F_{x_1x_1}(P) = 0$
   $\Longrightarrow$ $F_{x_0x_1}(P) = 0$, $F_{x_0x_0}(P) = 0$.}\]

(2) Let $a \ne b \in A$. 
As is well known of a property of vector spaces, 
\[\dim (\cH_a \cap \cH_b) 
   = \dim \cH_a + \dim \cH_b - \dim (\cH_a+\cH_b).\]
In the case of convex cones, 
\[\dim (\cP_a \cap \cP_b) 
   \leq \dim \cP_a + \dim \cP_b - (\cH_a+\cH_b)\]
is true, but $\leq$ cannot be replaced by $=$ in general. 
So, we must be careful to compute 
$\dim \big(\cP_1 \cap \cdots \cap \cP_r\big)$. 
\endproclaim

\def\Tqacnode{Definition 2.16} 
\proclaim{Definition 2.16} 
{\rm For an irreducible curve $C = V_{\C}(f)$ with $f \in \cH_{3,d}$, 
we say $C$ has a (simple) node at $P$ 
if two analytic branches of $C$ intersect at $P$ transversally. 
We say $C$ has an {\it acnode} at $P$, 
if $P$ is a simple node of $V_{\C}(f)$ and $P$ is an 
isolated point of $V_{\R}(f)$. 
In other words, an acnode is a real node whose tangents are 
non-real complex conjugates. }
\endproclaim

Note that if $f \in \cP(A$, $\cH)$, $f(P)=0$, $f$ is irreducible, 
and $V_{\C}(f)$ has a node at $P \in \Int(A)$, 
then $P$ is an acnode of $V_{\C}(f)$. 
Theorem 1.1 can be restated as the following using the notion of 
infinitely near zeros. 

\def\Tqbaf{Theorem 2.17} 
\proclaim{Theorem 2.17} 
{\sl Assume that $f \in \cP_{3,6}$ is not a square of a cubic polynomial. 
Then, $f \in \cE(\cP_{3,6})$ if and only if $f$ has just $10$ zeros 
on $\P_{\R}^2$ including all the infinitely near zeros.} 
\endproclaim

\Proof
Assume that $f \in \cE(\cP_{3,6})$ is not 
a square of cubic polynomial. 
Then $f$ is a limit of a sequence $\{f_n\}$ of exposed extremal elements 
in $\cP_{3,6}$ (see \cite{RefCLR}). 
Each $f_n$ has distinct $10$ zeros. 
Any infinitely near zero of $f$ is a limit of a squence of a certain 
zero of $f_n$. 
Thus, $f$ also has just $10$ zeros including all infinitely near zeros. 

Assume that $f$ has $10$ zeros $P_1$,$\ldots$, $P_{10}$ 
including infinitely near zeros. 
Consider the case that $f$ is irreducible in $\C[x_0,x_1,x_2]$. 
Let $\psi \colon X \to \P_{\C}^2$ be a proper birational morphism 
such that we can regard $P_1$,$\ldots$, $P_{10}$ are 
distinct points in $X$ in Noether's sense. 
As is well known in algebraic geometry, 
there exists a unique irreducible sextic curve 
$C \subset \P_{\C}^2$ such 
that the strict transform of $C$ to $X$ has 
nodes at $P_1$,$\ldots$, $P_{10}$. 
Thus $C = V_{\C}(f)$ and $f$ is extremal.  

Next consider the case that $f$ is reducible. 
Then $f = f_1 f_2$ ($f_1 \in \cP_{3,2}$, $f_2 \in \cP_{3,4}$) or 
$f = f_3 \overline{f_3}$ ($f_3$ is a cubic). 

If $f = f_1 f_2$, then $V_{\R}(f) = V_{\R}(f_1) \cup V_{\R}(f_2)$. 
If $V_{\R}(f_1)$ and $V_{\R}(f_2)$ are finite sets, 
then $\# V_{\R}(f_1) \leq 1$ and $\# V_{\R}(f_2) \leq 4$. 
Thus $f$ cannot have $10$ zeros. 

Consider the case $f = f_3 \overline{f_3}$. 
We assume that $V_{\R}(f_3)$ is a finite set. 
Then $f_3(P)=0$ for $P \in \R^3$ if and only if $\overline{f_3}(P)=0$. 
Thus $P \in V_{\R}(f_3) \cap V_{\R}(\overline{f_3}) 
  \subset V_{\C}(f_3) \cap V_{\C}(\overline{f_3})$. 
But the intersection of two distinct cubic curves $V_{\C}(f_3)$ 
and $V_{\C}(\overline{f_3})$ consists of $9$ points. 
Thus $f$ cannot have $10$ zeros. 
\end{proof}

Using the ideas in \cite{RefCLR, RefCLRb}, 
we also obtain the following theorem. 

\def\Tqbag{Theorem 2.18} 
\proclaim{Theorem 2.18}
{\sl Assume $d \geq 3$, and $f \in \cE(\cP_{3,2d})$ is irreducible. 
Let $N$ be the numbers of zeros of $f$ in $\P_{\R}^2$ including 
all the infinitely near zeros. Then 
\[\frac{(d+1)(d+2)}{2} \leq N \leq (2d-1)(d-1).\]

}
\endproclaim

\Proof
Let $P_1$,$\ldots$, $P_N$ be all the zeros of $f$ on $\P_{\R}^2$ 
including infinitely near zeros. 
There exists local cones or infinitesimal local cones 
$\cL_1$,$\ldots$, $\cL_r \subset \cP_{3,2d}$ 
such that $\cL_1 \cap \cdots \cap \cL_r = \R_+ \cdot f$. 
We may assume that $r=N$ and $\cL_i$ corresponds to $P_i$. 

If $N < (d+1)(d+2)/2$, then there exists $g \in \cH_{3,d}$ such that 
$P_1$,$\ldots$, $P_N \in V_{\R}(g)$ in the sense of Noether. 
Then $g^2 \in \cL_1 \cap \cdots \cap \cL_r = \R_+ \cdot f$. 
This implies $g^2 = cf$ ($\exists c \in \R_+$), and $f$ is reducible. 
Thus $N \geq (d+1)(d+2)/2$. 

Let $C := V_{\C}(f) \subset \P_{\C}^2$. Then 
$p_a(C) \geq \sum_{P} \nu(P)(\nu(P)-1)/2 + g(C') \geq N + g(C')$, 
where $p_a(C) = (2d-1)(2d-2)/2$ is the arithmetic (or virtual) genus, 
$g(C')$ is the genus of the normalization $C'$ of $C$, 
and $\nu(P_i)$ is the multiplicity of $C$ at $P_i$. 
Thus, we have $N \leq (2d-1)(d-1)$. 
\end{proof}

\section{Extremal elements of $\cP_{3,3}^+$}

In \S 3 and \S 4, 
we usually use the symbol $(x_0 \colon x_1 \colon x_2)$ to denote 
the standard homogeneous coordinate system of $\P_+^2$, $\P_{\R}^2$ 
or $\P_{\C}^2$. 
We sometime rewrite $(x_0 \colon x_1 \colon x_2)$ by $(x \colon y \colon z)$. 
But we also often use $(x$, $y)$ to denote a local coordinate system 
when there is no fear of confusion. 

\bigbreak
\subsection{Reducible elements of $\cE(\cP_{3,3}^+)$}\hfill\break\relax
\hbox{}\hskip\parindent 
Let $\XYZ := \big(\R_+ \cdot x \cup \R_+ \cdot y \cup \R_+ \cdot z\big) 
   - \{0\} \subset \cH_{3,1}$. 

\def\Ttaaa{Proposition 3.1}
\proclaim{Proposition 3.1} 
{\sl If $f \in \XYZ$, $g \in \cH_{3,d}$ 
and $fg \in \cE(\cP_{3,d+1}^+)$, 
then $g \in \cE(\cP_{3,d}^+)$. 
Conversely, if $f \in \XYZ$ and $g \in \cE(\cP_{3,d}^+)$, 
then $fg \in \cE(\cP_{3,d+1}^+)$. }
\endproclaim

\Proof
Assume that $fg \in \cE(\cP_{3,d+1}^+)$ and 
$g = h_1 + h_2$ for certain $h_1$, $h_2 \in \cP_{3,d}^+ - \{0\}$. 
Since $fg$ is extremal, there exists $c_1$, $c_2 \in \R_+$ such that 
$fh_1 = c_1fg$ and $fh_2=c_2fg$. Thus, $h_1 = c_1g$ and $h_2 = c_2g$. 
That is, $g$ is extremal. 

Conversely, we assume $f \in \XYZ$, $g \in \cE(\cP_{3,d}^+)$, 
and $fg = h_1 + h_2$ for certain $h_1$, $h_2 \in \cP_{3,d+1}^+ - \{0\}$. 
Note that $V_+(h_1+h_2) = V_+(h_1) \cap V_+(h_2)$, 
because $h_1 \geq 0$ and $h_2 \geq 0$ on $\P_+^2$. 
Thus $V_+(f) \subset V_+(h_1+h_2) \subset V_+(h_1) \cap V_+(h_2)$. 
Since $\dim V_+(f) = 1$, we have $V_{\C}(f) \subset V_{\C}(h_i)$ ($i=1$, $2$). 
Since $f$ is irreducible, $h_1$ and $h_2$ must be multiples of $f$. 
Let $h_1 = fg_1$ and $h_2 = fg_2$ ($g_1$, $g_2 \in \cH_{3,d}$). 
Since $h_1$, $h_2 \in \cP_{3,d+1}^+ - \{0\}$, 
we have $g_1$, $g_2 \in \cP_{3,d}^+ - \{0\}$. 
Since $g = g_1+g_2$ is extremal, 
there exists $c_1$, $c_2 \in \R_+$ such that $g_1 = c_1g$ and $g_2=c_2g$. 
Thus, $fg$ is extremal. 
\end{proof}

\def\Ttaba{Proposition 3.2}
\proclaim{Proposition 3.2} 
{\sl If $f \in \cE(\cP_{3,2}^+)$, 
then one of the following statements holds. } \par
{\parindent=20pt
\Item{(1)} {\sl $f = f_1 f_2$ where $f_1$, $f_2 \in \XYZ$.} 
\Item{(2)} {\sl There exists $f_1 \in \cH_{3,1}$ such that 
$f = f_1^2$ and $V_{\R}(f_1) \cap \Int(\P_+^2) \ne \emptyset$.}

}
\endproclaim

\Proof
Since $f$ is extremal, there exists $P \in \P_+^2$ such that $f(P) = 0$ 
by \Tqbb(1). 

(i) Assume that $P \in \Int(\P_+^2)$ and $V_+(f) = \{P\}$. 
By the classification of quadratic curves, 
this occurs only in the case $f = g_1^2 + g_2^2$ 
where $g_1$, $g_2 \in \cH_{3,1}$. 
In this case $V_{\R}(g_1)$ and $V_{\R}(g_2)$ are distinct lines 
which intersect at $P$. 
Thus, $f$ is not extremal. 

(ii) Assume that $P \in \Int(\P_+^2)$ and $V_+(f) \supsetneq \{P\}$. 
Then, there exists $Q \in \P_+^2$ such that $f(Q) = 0$ and $P \ne Q$. 
Since $f$ is PSD, $V_+(f)$ cannot be a real conic. 
Thus $V_+(f)$ must be a line passing through $P$ and $Q$. 
Thus (2) occurs. 

(iii) Assume that $V_{\R}(f) \cap \Int(\P_+^2) = \emptyset$ 
and $P \in \partial \P_+^2$. 
We may assume $P = (a \colon 0 \colon 1)$ where $a \geq 0$. 

It is easy to see that 
$\dim (\cP_{3,2}^+)_P \geq \dim \cP_{3,2}^+ - 2 = 4$. 
Since $f \in (\cP_{3,2}^+)_P$ is extremal, 
there exists $Q \in \partial \P_+^2$ such that $f(Q) = 0$. 

If $Q = (b \colon 0 \colon 1)$ ($a \ne b$, $b>0$), 
then, $f = c y^2$ ($\exists c > 0$). 
If $Q = (0 \colon b \colon 1)$ ($b>0$), 
then, $f = c x y$ ($\exists c > 0$). 
If $Q = (1 \colon b \colon 0)$ ($b>0$), 
then, $f = c y z$ ($\exists c > 0$). 
\end{proof}

\def\Ttaca{Proposition 3.3}
\proclaim{Proposition 3.3} 
{\sl Let $f \in \cE(\cP_{3,3}^+)$. 
If $f$ is reducible, then one the following statements holds. } \par
{\parindent=20pt
\Item{(1)} {\sl $f = f_1 f_2 f_3$ and $f_1$, $f_2$, $f_3 \in \XYZ$.} 
\Item{(2)} {\sl $f = f_1 f_2^2$, $f_1 \in \XYZ$, 
$f_2 \in \cH_{3,1}$  
and $V_{\R}(f_2) \cap \Int(\P_+^2) \ne \emptyset$.

}
{\sl Conversely, if $f \in \cH_{3,3}$ satisfies (1) or (2), 
then $f \in \cE(\cP_{3,3}^+)$. }
\endproclaim

\Proof
(i) Since $f$ is reducible, we can write as $f = f_1 g$, where 
$f_1 = ax+by+cz$ ($a$, $b$, $c \in \R$) and $g \in \cH_{3,2}$. 

(i-1) Consider the case $V_{\R}(f_1) \cap \Int(\P_+^2) = \emptyset$.
We may assume that $f_1 \geq 0$ on $\P_+^2$. 
Then $a \geq 0$, $b \geq 0$, $c \geq 0$ and $g \in \cP_{3,2}^+$. 
Assume that $a>0$ and $b>0$. 
Then $f = a x g + b y g + c z g$ and 
$a x g$, $b y g$, $c z g \in \cP_{3,3}^+$. 
Thus $f$ is not extremal. 
This implies $f_1 \in \XYZ$. 
We may assume $f_1 = x$. 
Since $f$ is extremal in $\cP_{3,3}^+$, 
$g$ must be extremal in $\cP_{3,2}^+$ by \Ttaaa. 
Then, we have the conclusion by \Ttaba.  

(i-2) Consider the case $V_{\R}(f_1) \cap \Int(\P_+^2) \ne \emptyset$.
Then $f$ must be divisible by $f_1^2$, because $f \geq 0$ on $\P_+^2$. 
So, $f = f_1^2 f_2$ where $f_2 \in \cH_{3,1}$ with 
$V_{\R}(f_2) \cap \Int(\P_+^2) = \emptyset$.
If we put $(f_2$, $f_1^2)$ as new $(f_1$, $g)$, 
we can apply (i-1). 

(ii) We prove the converse part. 
For $c>0$, the forms $cx^3$, $cx^2y$, $cxyz$ are also extremal 
in $\cP_{3,3}^+$. 
Thus if $f$ satisfies (1), then $f$ is extremal. 

Consider the case $f$ satisfies (2). 
Assume that $f = g_1+g_2$ where $g_1$, $g_2 \in \cP_{3,3}^+$. 
Then $V_+(f_1) \subset V_+(g_1) \cap V_+(g_2)$. 
Since $\dim V_+(f_1)=1$ and $f_1$ is irreducible, 
there exists $h_1$, $h_2 \in \cP_{3,2}^+$ such that 
$g_1 = f_1 h_1$, $g_2 = f_1 h_2$. 
Since $V_+(f_2) \subset V_+(h_1) \cap V_+(h_2)$, 
there exits $c_1$, $c_2 \in \R_+$ such that $h_1 = c_1 f_2^2$, 
$h_2 = c_2 f_2^2$.  Thus $f \in \cE(\cP_{3,3}^+)$. 
\end{proof}

\def\Ttada{Proposition 3.4}
\proclaim{Proposition 3.4} 
{\sl If $f \in \cH_{3,3}$ is 
irreducible in $\R[x,y,z]$, then $f$ is irreducible in $\C[x,y,z]$. }
\endproclaim

\Proof
Assume that $f$ is divisible by $f_1 \in \C[x,y,z]$ with 
$f_1 \notin \R[x,y,z]$. Let $\overline{f_1}$ be the complex conjugate 
of $f_1$. Then, $f$ is divisible by $\overline{f_1}$. 
Thus $f = f_1 \overline{f_1} f_2$. 
Since $f_1 \overline{f_1} \in \R[x,y,z]$, 
we have $f_2 \in \R[x,y,z]$. 
\end{proof}

\bigbreak
\subsection{Basic lemmata for irreducible elements of $\cE(\cP_{3,3}^+)$}\hfill\break\relax
\hbox{}\hskip\parindent 
In \S 3.2 --- \S 4, we use the symbols 
$P_x = P_{x_0} := (1 \colon 0 \colon 0)$, 
$P_y = P_{x_1} := (0 \colon 1 \colon 0)$, 
$P_z = P_{x_2} := (0 \colon 0 \colon 1) \in \P_+^2$, 
$L_x = L_{x_0} := V_+(x) - \{P_y$, $P_z\}$, 
$L_y = L_{x_1} := V_+(y) - \{P_z$, $P_x\}$, 
and $L_z = L_{x_2} := V_+(z) - \{P_x$, $P_y\}$. 
But meanings of symbols $Q_x$, $Q_y$ and $Q_z$ change 
in respective thorems, propositions and lemmata. 

For an irreducible complex algebraic curve $C$, we denote
\begin{align*}
 & \Sing(C) := \{P \in C \; \big| \; 
      \hbox{$P$ is a singular point of $C$}\big\}, \\
 & \Reg(C) := \{P \in C \; \big| \; 
      \hbox{$P$ is a non-singular point of $C$}\big\}. 
\end{align*}

\def\Ttaea{Lemma 3.5}
\def\Ttaean{3.5}
\proclaim{Lemma 3.5} 
{\sl Assume that $f \in \cE(\cP_{3,3}^+)$ is irreducible. 
Then, $V_{\C}(f)$ is a rational curve on $\P_{\C}^2$ whose unique 
singular point lies on $\P_+^2$. }
\endproclaim

\Proof
Assume that $\Sing(V_{\C}(f)) \cap \P_+^2 = \emptyset$. 
Then $f(P) > 0$ for all $P \in \Int(\P_+^2)$. Let 
\[D_i := \big\{ (x_0 \colon x_1 \colon x_2) \in \P_+^2 \; \big| \; 
  \hbox{$x_0^2+x_1^2+x_2^2 \leq 5x_i^2$} \big\}\]
($i=0$, $1$, $2$). 
Note that $\P_+^2 = D_0 \cup D_1 \cup D_2$. 
On $D_0$, we put $x := x_1/x_0$, $y := x_2/x_0$, 
and $f_0(x,y) := f(1,x,y)$. Let 
\[D := \big\{ (x,y) \in \R^2 \; \big| \; 
  \hbox{$x^2+y^2 \leq 4$, $x \geq y \geq 0$} \big\}.\]

(1) We shall prove that there exists $c_0 > 0$ 
such that $f_0(x,y) \geq c_0xy$ for all $(x,y) \in D$. 

Note that if $f_0(a,0) = 0$ for a certain $a > 0$, 
then $x$-axis is a tangent to $V_{\C}(f_0)$ at $(a$, $0)$. 
Since any cubic curve has no bitangent line, 
there exists at most one $a_0 > 0$ such that $f_0(a_0$, $0) = 0$. 

(1-1) Assume that $f_0(x,0) > 0$ for all $0 \leq x \leq 2$. 
Then there exists no $(x$, $y) \in D$ such that $f_0(x$, $y) = 0$. 
Thus 
$m := \min \big \{ f_0(x,y)$ $\big|$ $(x,y) \in D \big\} > 0$. 
Put $c_0 := m/4$, then $f_0(x,y) \geq c_0xy$ for all $(x,y) \in D$. 

(1-2) Assume that $f_0(a,0) = 0$ for a certain $0 \leq a \leq 2$. 
There exist $c_1$,$\ldots$, $c_8 \in \R$ such that 
\[f_0(x,y) = c_1y + c_2(x-a)^2 + 2c_3(x-a)y + c_4y^2 
   + g(x-a,y)\]
where $g(s,t) = c_5s^3 + c_6s^2t + c_7st^2 + c_8t^3$. 
Since $f_0(a,y) \geq 0$ for all $y \geq 0$, we have $c_1 \geq 0$. 
If $c_1 = 0$, then $(a$, $0)$ is a singular point of $V_{\C}(f)$. 
Thus $c_1 > 0$. 
Then, there exists an open neighborhood $(a$, $0) \in U_a \subset D$ such 
that $f_0(x,y) \geq (c_1/4) y$ for all $(x$, $y) \in U_a$. 
Then $f_0(x,y) \geq (c_1/8) xy$ for all $(x$, $y) \in U_a$. 
So, we put $m_1(a) := c_1/8$. 

Let $V := \big\{a \in [0$, $2]$ $\big|$ $f_0(a$, $0) = 0\big\}
  \subset \{0$, $a_0\}$, 
$\displaystyle U := \bigcup_{a \in V} U_a$, and 
$\displaystyle m_2 := \min_{a \in V} m_1(a)$. 
Then $f_0(x$, $y) \geq m_2xy$ for all $(x$, $y) \in U$. 
Note that
\[m_3 := \min \big\{ f_0(x,y) \; \big| \; (x,y) \in \Cls(D-U) \big\} > 0.\]
So, put $c_0 := \min\{m_2$, $m_3/4\}$. 
Then $f_0(x,y) \geq c_0xy$ for all $(x,y) \in D$. 

\smallskip

By (1), there exists $c > 0$ such that $F(x_0,x_1,x_2) \geq cx_0x_1x_2$ for 
all $(x_0 \colon x_1 \colon x_2) \allowbreak \in \P_+^2$. 
So, $f(x_0,x_1,x_2) - cx_0x_1x_2 \in \cP_{3,3}^+$ 
and $f$ is not extremal. 
A contradiction. 
Thus $\Sing(V_{\C}(f)) \cap \P_+^2 \ne \emptyset$. 

Since $V_{\C}(f)$ is a cubic curve, $p_a(C) = 1$. 
Since $\Sing(V_{\C}(f)) \cap \P_+^2 \ne \emptyset$, 
$V_{\C}(f)$ has just one singular point $P$ and $V_{\C}(f)$ is 
a rational curve, by Riemann  genus formula. 
Thus, $P \in \P_+^2$. 
\end{proof}

\def\Ttafb{Lemma 3.6}
\proclaim{Lemma 3.6} 
{\sl Assume that $f \in \cE(\cP_{3,3}^+)$ is irreducible, 
and $P$ is the unique singular point of $V_{\C}(f)$. 
Then $P \notin \big\{P_x$, $P_y$, $P_z \big\}$. }
\endproclaim

\Proof
Assume that $P \in \big\{P_x$, $P_y$, $P_z \big\}$. 
We may assume $P = P_x$. 
We use the same notation as in the proof of the previous lemma. 
Then $P = (0$, $0) \in D \subset D_0 \subset \P_+^2$. 
Then $f_0(x,y) = g_2(x,y) + g_3(x,y)$ where 
$g_d(x,y)$ is a homogeneous polynomial of degree $d$. 
Then $f(x_0,x_1,x_2) = g_2(x_1,x_2)x_0 + g_3(x_1,x_2)$. 
Considering the cases $x_0=0$ and $x_0 \to +\infty$, we conclude that 
$g_2(x_1,x_2) \in \cP_{2,2}^+$ and $g_3(x_1,x_2) \in \cP_{2,3}^+$. 
If $g_3=0$, then $f=x_0g_2$ is not irreducible. So, $g_3 \ne 0$. 
If $g_2=0$, the cubic $g_3(x_1$, $x_2)$ is reducible in $\C[x_1$, $x_2]$. 
Thus, $g_3 \ne 0$. 
Therefore, $f$ is not extremal. 
\end{proof}

\def\Ttagb{Lemma 3.7}
\proclaim{Lemma 3.7} 
{\sl Assume that $f \in \cE(\cP_{3,3}^+)$ is irreducible. 
Then $I_{P_y}(f$, $x) + I_{P_z}(f$, $x) \leq 3$. 
If $V_+(f) \cap L_x \ne \emptyset$, then $\# (V_+(f) \cap L_x) = 1$ and 
$I_{P_y}(f$, $x) + I_{P_z}(f$, $x) \leq 1$. }
\endproclaim

\Proof
Assume that $f \in \cE(\cP_{3,3}^+)$ is irreducible. 
Then $I_{P_y}(f$, $x) + I_{P_z}(f$, $x) \leq I(f,x) = \deg f = 3$. 
If $Q \in V_+(f) \cap L_x \ne \emptyset$, 
then, $I_{P_y}(f$, $x) + I_{P_z}(f$, $x) + I_Q(f$, $x) \leq I(f,x) = 3$. 
Since $I_Q(f$, $x) \geq 2$, we have 
$I_{P_y}(f$, $x) + I_{P_z}(f$, $x) \leq 1$. 
\end{proof}

\def\Ttag{Lemma 3.8}
\proclaim{Lemma 3.8} 
{\sl Let $\cP = \cP_{3,3}^+$. 
Assume that $f(x$, $y$, $z) \in \cE(\cP)$ is 
irreducible in $\C[x$, $y$, $z]$. 
Then, $V_+(f)$ is a finite set, 
and we can choose local cones or infinitesimal local cones $\cL_i$ 
of $\cP$ at $P_i$ which satisfy the following conditions. }\par
{\parindent=20pt
\Item{(1)} {\sl $\cL_1 \cap \cdots \cap \cL_r = \R_+ \cdot f$. }
\Item{(2)} {\sl $P_i \ne P_j$ if $i \ne j$. }
\Item{(3)} {\sl $\{P_1$,$\ldots$, $P_r\} = V_+(f)$. }
\Item{(4)} {\sl If $P_i \in V_+(f) \cap \Reg(V_{\C}(f))$, 
then $P_i \in \partial \P_+^2$, 
and one of the following cases occurs. 
Here $x_0=x$, $x_1=y$ and $x_2=z$.}
{\parindent=40pt
\Item{(4-1)} {\sl If $P \in L_{x_j}$ and 
$k \in \{0$, $1$, $2\} - \{j\}$, then }
\[\cL_i = \big\{ F \in \cP \; \big| \; \hbox{$F(P)=F_{x_k}(P)=0$} \big\}.\]
\Item{(4-2)} {\sl Assume that $P=P_{x_j}$, 
$\{j$, $k$, $l\} = \{0$, $1$, $2\}$ 
and $m := I_P(f$, $x_k) \geq I_P(f$, $x_l) = 1$. 
If $m=1$, then 
\[\cL_i = \big\{ F \in \cP \; \big| \; \hbox{$F(P)=0$} \big\}.\]
If $m=2$, then 
\[\cL_i 
 = \big\{ F \in \cP \; \big| \; \hbox{$F(P)=F_{x_k}(P)=0$} \big\}.\]
If $m=3$, then }
\[\cL_i = \big\{ F(x,y,z) \in \cP \; \big| \; 
      \hbox{$F(P)=F_{x_k}(P)=F_{x_kx_k}(P)=0$} \big\}.\]

}
\Item{(5)} {\sl If $P_i \in V_+(f)$ is an acnode, then }
\[\cL_i 
 = \big\{ F(x,y,z) \in \cP \; \big| \; 
    \hbox{$F(P)=F_x(P)=F_y(P)=0$} \big\}.\]

}
\endproclaim

\Proof
(1) By \Tqbae, 
there exist local cones or infinitesimal local cones $\cL_i$ 
of $\cP$ at $P_i$ which satisfy 
$\cL_1 \cap \cdots \cap \cL_r = \R_+ \cdot f$. 
If $\cL_i \subset \cL_j$ for $i \ne j$, we get rid of $\cL_j$. 
So, we may assume that $\cL_i \not\subset \cL_j$ if $i \ne j$. 

\smallskip

(2) By our observation in \Tqbi, the local cone $(\cP_{3,3}^+)_{P_i}$ 
does not have two distinct infinitesimal local cones, 
if $f$ is cubic and $P$ is an acnode, a cusp or a non-singular point in 
$\partial \P_+^2$. Thus, $P_i \ne P_j$ if $i \ne j$. 

\smallskip

(3) $P_i \in V_+(f)$ by the definition. 
Assume that $P_{r+1} \in V_+(f) - \{P_1$,$\ldots$, $P_r\}$. 
Let $\cL_{P_{r+1}}$ be the local cone or the infinitesimal local cone 
of $\cP_{3,3}^+$ at $P_{r+1}$ with respect to $f$. 
$\cL_{P_{r+1}}$ is unique by \Tqbi. 
Then $\cL_{P_1} \cap \cdots \cap \cL_{P_{r+1}} = \R_+ \cdot f$ 
still holds. 
After repeating this process, we may assume 
$V_+(f) = \{P_1$,$\ldots$, $P_{r+k}\}$. 

\smallskip

(4) Assume that $P \in V_+(f) \cap \Reg(V_{\C}(f))$. 
If $P \in \Int(\P_+^2)$, then $V_+(f)$ must contain a locus of 
analytic curve $C_P$ near $P$, 
and $\sign(f)$ is opposite on both sides of $C_P$. 
Thus $P \in \partial \P_+^2$. 

\smallskip

(4-1) Consider the case $P \in L_z$. 
Note that $I_P(f$, $z) \leq I(f$, $z) = 3$. 
If $I_P(f$, $z)$ is odd, then $V_+(f)$ must contain a locus of 
analytic curve $C_P$ near $P$, and $\sign(f)$ changes across $C_P$. 
Thus $I_P(f$, $z) = 2$. 
By \Tqbi(3), we have (4-1). 

\smallskip

(4-2) follows from \Tqbi(4). 

\smallskip

(5) follows from \Tqbi(2). 
\end{proof}

\def\Ttagc{Notation 3.9}
\proclaim{Notation 3.9} 
{\rm Throughout the remaining part of \S 3, we use the following notation. 
We denote a local cone or an infinitesimal local cone $\cL_i$ 
in the above lemma by $\cL_{P_i}$. 
Assume that $f(x$, $y$, $z) \in \cE(\cP)$ is irreducible in $\C[x$, $y$, $z]$. 
\Ttag(1) is represented as 
\[\cL_{P_1} \cap \cdots \cap \cL_{P_r} = \R_+ \cdot f\]
where $\{P_1$,$\ldots$, $P_r \} = V_+(f)$ and $P_i \ne P_j$ if $i \ne j$. 
This $\cL_{P_1} \cap \cdots \cap \cL_{P_r}$ is denoted by $\cL_f$. }
\endproclaim

\def\Ttaic{Lemma 3.10}
\def\Ttaicn{3.10}
\proclaim{Lemma 3.10} 
{\sl Assume that $f \in \cE(\cP_{3,3}^+)$ is irreducible. 
Let $\pi \colon \P_{\R}^2 \to \P_+^2$ be the surjective morphism defined by 
$\pi(x_0 \colon x_1 \colon x_2) = (x_0^2 \colon x_1^2 \colon x_2^2)$, 
and let $g := \pi^* f$. That is $g(x_0,x_1,x_2) = f(x_0^2$, $x_1^2$, $x_2^2)$. 
Take $P \in V_+(f)$. We define 
\[N_P(f) := \left\{
\begin{array}{ll}
  4 & \hbox{(If $P$ is an acnode of $V_{\C}(f)$.)} \\
 12 & 
   \hbox{(If $P \in \Int(\P_+^2)$ and $P$ is a cusp of $V_{\C}(f)$.)} \\
  6 & 
  \hbox{(If $P \in \partial \P_+^2$ and $P$ is a cusp of $V_{\C}(f)$.)} \\
  2 & \hbox{(If $P \in L_{x_0} \cup L_{x_1} \cup L_{x_2}$ 
     and $P \in \Reg(V_{\C}(f))$.)} \\
  \length_P f & 
     \hbox{(If $P \in \{P_{x_0}$, $P_{x_1}$, $P_{x_2}\}$ 
     and $P \in \Reg(V_{\C}(f))$.)}
\end{array}\right.\]
and $\displaystyle N(f) := \sum_{Q \in V_+(f)} N_Q(f)$. 
Then $N(f) \leq 10$. 
Moreover, if $N(f)=10$, then $g$ is irreducible in $\C[x_0,x_1,x_2]$ 
and $g \in \cE(\cP_{3,6})$. 
In particular, $V_{\C}(f)$ does not have a cusp in $\Int(\P_+^2)$. }
\endproclaim

\Proof
We will separate 8 cases according to the type of the point $P$. 

(1) Consider the case $P \in \Int(\P_+^2)$ and $P$ is 
the acnode of $V_{\C}(f)$. 
If $P = (1 \colon p^2 \colon q^2)$ ($p>0$, $q>0$), 
then $\pi^{-1}(P) = \{(1 \colon p \colon q)$, $(1 \colon p \colon -q)$, 
$(1 \colon -p \colon q)$, $(1 \colon -p \colon -q)\}$. 
Put $x := (x_1 - p^2 x_0)/x_0$ and $y := (x_2 - q^2 x_0)/x_0$. 
Then $f(x,y) = ax^2 + 2bxy + cy^2 + (\hbox{\rm higher terms})$ as in \Tqbi(1). 
Put $x' := (x_1 - p x_0)/x_0$ and $y' := (x_2 - q x_0)/x_0$ 
around $Q := (1 \colon p \colon q) \in \P_{\R}^2$. 
Then $g(x',y') = ax'^2 + 2bx'y' + cy'^2 + (\hbox{\rm higher terms})$. 
Thus, $\length_Q g = 1$ for all $Q \in \pi^{-1}(P)$. 
Therefore, $V_{\C}(g)$ has $4$ acnodes in $\pi^{-1}(P)$. 

(2) Consider the case $P \in L_x \cup L_y \cup L_z$ is 
the acnode of $V_{\C}(f)$. 
We may assume that $P = (1 \colon p^2 \colon 0)$ ($p>0$). 
Let $x := (x_1 - p^2 x_0)/x_0$, $y := x_2/x_0$, 
$x' := (x_1 - p x_0)/x_0$ and $y' := x_2/x_0$. 
Take $Q := (1 \colon p \colon 0) \in \pi^{-1}(P) 
= \{(1 \colon \pm p \colon 0)\}$. 
Then $g(x',y') = ax'^2 + 2bx'y'^2 + cy'^4 + (\hbox{\rm higher terms})$. 
Thus $\length_Q g = 2$. 
Therefore, $V_{\C}(g)$ has $4$ zeros or 
infinitely near zeros in $\pi^{-1}(P)$. 

(3) Consider the case $P \in \Int(\P_+^2)$ and $P$ is the cusp of $V_{\C}(f)$. 
If $P = (1 \colon p^2 \colon q^2)$ ($p>0$, $q>0$), 
then $f(x,y) = x^3 + y^2 + (\hbox{\rm higher terms})$ as in \Tqbi(5). 
Then $g(x',y') = x'^3 + y'^2 + (\hbox{\rm higher terms})$. 
Thus $\length_Q g = 3$. 
Therefore, $V_{\C}(g)$ has $12$ zeros or 
infinitely near zeros in $\pi^{-1}(P)$. 

\smallskip

(4) Consider the case $P \in L_x \cup L_y \cup L_z$ 
and $P$ is the cusp of $V_{\C}(f)$. 
Assume that $P = (1 \colon p^2 \colon 0)$ ($p>0$), 
Then $g(x',y') = x'^6 + y'^2 + (\hbox{\rm higher terms})$. 
Thus $\length_Q g = 3$. 
Therefore, $V_{\C}(g)$ has $6$ zeros or 
infinitely near zeros in $\pi^{-1}(P)$. 

\smallskip

(5) Consider the case $P \in L_x \cup L_y \cup L_z$ and 
$P \in \Reg(V_{\C}(f))$. 
Then $f(x,y) = a y + b x^2 + (\hbox{\rm higher terms})$ ($a>0$, $b>0$) 
as in \Tqbi(3). 
Since $g(x',y') = a y'^2 + b x'^2 + (\hbox{\rm higher terms})$, 
we have $\length_Q g = 1$. 
Thus, $V_{\C}(g)$ has $2$ acnodes in $\pi^{-1}(P)$. 

\smallskip

(6) Consider the case $P \in \{P_{x_0}$, $P_{x_1}$, $P_{x_2}\}$ 
and $P \in \Reg(V_{\C}(f))$. 
Let $x := x_1/x_0$ and $y := x_2/x_0$. 
We may assume that 
$f(x$, $y) = ax + by^m + (\hbox{\rm higher terms})$ with $a>0$, $b>0$. 
Then 
$g(x$, $y) = ax^2 + by^{2m} + (\hbox{\rm higher terms})$. 
Thus, $\length_Q g = \length_P f$. 

\smallskip

(7) Consider the case that $g = g_1g_2$ is reducible in $\C[x,y,z]$. 

Then $\pi(V_{\C}(g_1)) \cup \pi(V_{\C}(g_2)) = \pi(V_{\C}(g)) = V_{\C}(f)$. 
Since $V_{\C}(f)$ is an irreducible cubic curve, 
$V_{\C}(g_1)$ and $V_{\C}(g_2)$ are irreducible cubic curves. 
Thus, we may assume that 
there exists $g_3 \in \C[x,y,z]_{(3)} - \R[x,y,z]$ such 
that $g_1 = g_3$ and $g_2 = \overline{g_3}$. 
$g_3$ is irreducible in $\C[x,y,z]$, because $V_{\C}(g_3)$ is irreducible. 

The intersection number of $V_{\C}(g_3)$ and $V_{\C}(\overline{g_3})$ is 
equal to $9$. 
So, $N(f) \leq 9$. 

\smallskip

(8) Consider the case $N(f) \geq 10$. 

Then, $g$ is irreducible, and the arithmetic genus of $V_{\C}(g)$ satisfies 
\[10 = \frac{(\deg g - 1)(\deg g - 2)}{2}
   \geq \sum_{i=1}^n \frac{\nu(Q_i')(\nu(Q_i')-1)}{2} \geq n.\]
So, $N(f) = 10$. 
Then $g \in \cE(\cP_{3,6})$ by \Tqbaf. 
\end{proof}

\def\Ttah{Lemma 3.11}
\proclaim{Lemma 3.11} 
{\sl Assume that $f \in \cE(\cP_{3,3}^+)$ is irreducible, 
and $P$ is the acnode of $V_{\C}(f)$. 
Then $P \notin \partial \P_+^2$.} 
\endproclaim

\Proof
Assume that $P$ is an acnode of $V_{\C}(f)$ and $P \in \partial \P_+^2$. 
Then, $P \in L_x \cup L_y \cup L_z$ by \Ttaea \ and \Ttafbn. 
We may assume $P = (0 \colon a \colon 1) \in L_x$ ($a>0$). 
If we replace $y$ by $y/a$, we may assume $P = (0 \colon 1 \colon 1)$. 
$f(P_y)=f(P_z)=0$ is impossible by \Ttagb. 
So, we may further assume $f(P_y)>0$ by the symmetry. 

\smallskip

(1) Consider the case $Q_z = (1 \colon r \colon 0) \in V_+(f)$ ($r>0$). 

Then $V_+(f) \subset \{P$, $Q_z\} \cup V_+(y)$ by \Ttagb. 
In \Ttagc, we put $P_1=P$ and $P_2=Q_z$. 
Let $G(x$, $y$, $z) := y(r x-y-z)^2$. 
Then, $G_x=2r y(r x-y-z)$, $G_{xx}=2r^2y$, $G_y=(r x-y-z)(r x-3y-z)$, 
$G_z=-2y(r x-y-z)$ and $G_{zz}=2y$. 

For $Q \in \{P$, $Q_z\} = \{P_1$, $P_2\}$, 
we have $G(Q)=G_x(Q)=G_y(Q)=G_z(Q)=0$. 
Thus, $G \in \cL_{P_1} \cap \cL_{P_2}$. 

For any point $Q \in V_+(y)$, 
we have $G(Q)=G_x(Q)=G_{xx}(Q)=G_z(Q)=G_{zz}(Q)=0$. 
Note that $I_{P_x}(f$, $z) \leq 1$ by \Ttagb. 
Thus, if $P_x \in V_+(f)$, $\cL_{P_x}$ is one of the following: 
\begin{align*}
 & \cL_1 := \big\{ g \in \cP_{3,3}^+ \; \big| \; 
   \hbox{$g(P_x) = 0$} \big\}, \\
 & \cL_2 := \big\{ g \in \cP_{3,3}^+ \; \big| \; 
   \hbox{$g(P_x) = g_z(P_x) = 0$} \big\}, \\
 & \cL_3 := \big\{ g \in \cP_{3,3}^+ \; \big| \; 
   \hbox{$g(P_x) = g_z(P_x) = g_{zz}(P_x) = 0$} \big\}. 
\end{align*}
Since $G(P_x)=G_z(P_x)=G_{zz}(P_x)=0$, 
we have $G \in \cL_i$ ($i=1$, $2$, $3$). 
Thus, $G \in \cL_{P_x}$ if $P_x \in V_+(f)$. 

Similarly, since $I_{P_z}(f$, $x) \leq 1$, 
we have $G \in \cL_{P_z}$ if $P_z \in V_+(f)$. 

If $Q \in L_y \cap V_+(f)$, then 
\[\cL_Q = \big\{ g \in \cP_{3,3}^+ \; \big| \; 
  \hbox{$g(Q) = g_z(Q) = 0$} \big\}.\]
Thus, $G \in \cL_Q$. 

Since $V_+(f) \subset \{P$, $Q_z\} \cup V_+(y)$, we have 
$G \in \cL_{P_1} \cap \cdots \cap \cL_{P_r} = \R_+ \cdot f$. 
This implies $f$ is a multiple of $G$ and $f$ is reducible. 
A contradiction. 

\smallskip

(2) Consider the case $V_+(f) \cap L_z = \emptyset$. 

Then $V_+(f) \subset \{P\} \cup V_+(y)$. 
Let $H(x$, $y$, $z) := y(y-z)^2$. 
Then, $H_x=0$, $H_y=(y-z)(3y-z)$, $H_{yy}=2(3y-2z)$, 
$H_z=-2y(y-z)$ and $H_{zz}=2y$. 
Thus, for any point $Q \in V_+(y)$, 
we have $H(Q)=H_x(Q)=H_{xx}(Q)=H_z(Q)=H_{zz}(Q)=0$. 
Thus, if $Q \in V_+(f) \cap (L_y \cup \{P_z\})$, then $H \in \cL_Q$. 

Moreover, $H_y(P_x)=H_{yy}(P_x)=0$. 
So, $H \in \cL_{P_x}$. 
Therefore 
$H \in \cL_{P_1} \cap \cdots \cap \cL_{P_r} = \R_+ \cdot f$. 
A contradiction. 
\end{proof}

By \Ttaea --- \Ttahn, 
if $f \in \cE(\cP_{3,3}^+)$ is irreducible, 
$V_{\C}(f)$ has an acnode in $\Int(\P_+^2)$ or 
a cusp on $\partial \P_+^2 - \{P_x$, $P_y$, $P_z \}$. 

\bigbreak
\subsection{The case $V_{\C}(f)$ has a cusp}\hfill\break\relax
\hbox{}\hskip\parindent 
From this subsection, many complicated calculations appear.  
In most of them, we use the software Mathematica. 
The code for Mathematica can be found on the link of the authors WEB 
or in arXiv's anc folder. 

\def\Ttca{Lemma 3.12}
\proclaim{Lemma 3.12} 
{\sl Let $p \geq 0$, $q \geq 0$ be constants with $(p$, $q) \ne (0$, $0)$, 
and $P := (0 \colon 1 \colon 1) \in L_x \subset \P_+^2$. Let 
\[\cH_c := \left\{ g(x,y,z) \in \cH_{3,3} \ \left| \   
  \vcenter{
     \hbox{$g(P)=g_x(P)=g_y(P)=g_{xx}(P)=g_{xy}(P)=0,$}
     \hbox{$g(p,0,1)=g_x(p,0,1)=g(q,1,0)=g_y(q,1,0)=0$}}\right.\right\}.\]
Then $\cH_c = \R \cdot \frh_{pq}$, where 
\begin{align*}
 \frh_{pq}(x,y,z) 
 & := 2x^3 + 3(p-q) x^2 y
          - 6p q x y^2 + q^2(3p+q)y^3 
    + 3(q-p)x^2z - 6p q x z^2 \\
 & \hskip20pt 
    + (p^3+36p^2q-6p q^2-2q^3)y^2 z 
    + (-2p^3-6p^2q+3p q^2+q^3)y z^2 \\
 & \hskip20pt 
    + p^2(p+3q) z^3 + 12p q x y z. 
\end{align*}
Moreover, $\frh_{pq} \in \cE(\cP_{3,3}^+)$ 
and $V_{\C}(\frh_{pq})$ is irreducible. 
Moreover, $\frh_{pq}(x^2$, $y^2$, $z^2) \in \cE(\cP_{3,6})$. 
}
\endproclaim

\Proof
Let $\{e_1$,$\ldots$, $e_{10}\} = \{x^3$,$\ldots$, $xyz\}$ be 
all the monomials of $\cH_{3,3}$. 
Take $g = a_1e_1 + \cdots + a_{10} e_{10} \in \cH_{3,3}$, 
and let ${\bf a}$ be the column vector ${}^t(a_1$,$\ldots$, $a_{10})$. 
A differential equation $g(0$, $1$, $1)=g_x(0$, $1$, $1)
 = g_y(0$, $1$, $1)=g_{xx}(0$, $1$, $1)=g_{xy}(0$, $1$, $1)
 = g(p$, $0$, $1)=g_x(p$, $0$, $1)=g(q$, $1$, $0)=g_y(q$, $1$, $0) = 0$ 
can be written as $A{\bf a} = {\bf 0}$ for a certain $9 \times 10$ matrix $A$. 
Using Mathematica, we can check that 
$\Ker A = \R \cdot \frh_{pq}$. 
If $\frh_{pq}$ is PSD on $\P_+^2$, 
then $\frh_{pq} \in \cE(\cP_{3,3}^+)$ by \Tqbh. 

\smallskip

(1) Solving the equation $\frh_{pq}(x(t)$, $1$, $1+tx(t)) = 0$, 
then we have 
\[x(t) = -\frac{t^2(p+q)^3}{(p t-1)^2 ((p+3q)t+2)}.\]
Put $\displaystyle z(t) := 1 + t x(t) 
= -\frac{(q t+1)^2 ((3p+q)t-2)}{(p t-1)^2((p+3q)t+2)}$. 
Then $\frh_{pq}(x(t)$, $1$, $z(t))=0$ for all $t \in \C$. 
Thus $V_{\C}(\frh_{pq})$ is an irreducible rational curve 
which has a cusp at $P$. 

\smallskip

(2) We shall show that 
$V_{\R}(\frh_{pq}) \cap \Int(\P_+^2) = \emptyset$. 
It is enough to show that if $x(t)>0$ then $z(t) \leq 0$. 

If $x(t)>0$, then $(p+3q)t+2 < 0$ and 
$\displaystyle t \leq -\frac{2}{p+3q}$. Then 
\[(3p+q)t-2 \leq -(3p+q)\frac{2}{p+3q}-2 
    = -\frac{8(p+q)}{p+3q} < 0.\]
Thus $z(t) \leq 0$. 

\smallskip

(3) Since $\frh_{pq}(1$, $1$, $1) = 2p^3q^3 > 0$, 
we have $\frh_{pq} \in \cP_{3,3}^+$. 
Put $P_2 = (p \colon 0 \colon 1)$ and $P_3 = (q \colon 1 \colon 0)$. 
Then $N_P(\frh_{pq}) + N_{P_2}(\frh_{pq}) + N_{P_3}(\frh_{pq}) 
= 6+4+4 = 10$. 
Thus $\frh_{pq}(x^2$, $y^2$, $z^2) \in \cE(\cP_{3,6})$ 
by \Ttaic. 
\end{proof}

\def\Ttcb{Theorem 3.13}
\def\Ttcbn{3.13}
\proclaim{Theorem 3.13} 
{\sl Assume that $f \in \cE(\cP_{3,3}^+)$ is irreducible, 
and $V_{\C}(f)$ has a cusp at $P \in \P_{\C}^2$. 
Then $P \in L_x \cup L_y \cup L_z$. 

Assume that $P = (0 \colon 1 \colon 1)$. 
Then, there exists $p \geq 0$, $q \geq 0$ and $c>0$ such that 
$(p$, $q) \ne (0$, $0)$ and $f = c \frh_{pq}$. }
\endproclaim

\Proof
Let $P$ be a cusp of $V_{\C}(f)$. 
Then $P \in L_x \cup L_y \cup L_z$, by Lemmas \Ttaean, \Ttafbn \ and \Ttaicn. 
Assume that $P = (0 \colon 1 \colon 1)$. 
By \Tqbi(5), 
\[\cL_P 
   = \big\{ g(x,y,z) \in \cP_{3,3}^+ \; \big| \;  
       \hbox{$g(P)=g_x(P)=g_y(P)=g_{xx}(P)=g_{xy}(P)=0$}\big\}.\]
We consider conditions : \par
{\parindent=30pt
\Item{($Cy$)} $\exists P_2 = (p \colon 0 \colon 1) \in V_+(f)$ or 
``$I_{P_2}(f$, $y) \geq 2$ for $P_2 = (0 \colon 0 \colon 1)$ ($p:=0$)''. 
\Item{($Cz$)} $\exists P_3 = (q \colon 1 \colon 0) \in L_z \cap V_+(f)$ or 
``$I_{P_3}(f$, $y) \geq 2$ for $P_3 = (0 \colon 1 \colon 0)$ ($q:=0$)''. 

}
\smallskip

(1) Consider the case ($Cy$) and ($Cz$) are true. 
Then, $f = c \frh_{pq}$ ($p \geq 0$, $q \geq 0$, $c>0$) 
by \Ttca. 

\smallskip

(2) Consider the case ($Cy$) is true, ($Cz$) is false and $f(P_y) > 0$. 

Let $G(x,y,z) : = y(y-z)^2$. 
Then, $G_x = 0$, $G_{xx} = 0$, $G_{xy}=0$, 
$G_z = -2y(y-z)$ and  $G_{zz} = 2y$. 
Thus $G \in \cL_P$. 
Since $G(P_2)=G_x(P_2)=0$, we have $G \in \cL_{P_2}$. 

Note that $V_+(f) \subset \{P$, $P_2$, $P_x$, $P_z\}$ and 
$I_{P_z}(f$, $x) \leq 1$. 
Since $G(P_z)=G_x(P_z)=G_{xx}(P_z)=0$, we have $G \in \cL_{P_z}$. 
Since $G(P_x)=G_z(P_x)=G_{zz}(P_x)=0$, we have $G \in \cL_{P_x}$. 
Thus $G \in \cL_f = \R_+ \cdot f$. A contradiction. 

\smallskip

(3) Consider the case ($Cy$) is true, ($Cz$) is false and $f(P_y) = 0$. 
Let 
\[\cH_1 := \left\{ g(x,y,z) \in \cH_{3,3} \ \left| \  
   \vcenter{
      \hbox{$g(P)=g_x(P)=g_y(P)=g_{xx}(P)=g_{xy}(P)=0$,}
      \hbox{$g(p,0,1)=g_x(p,0,1)=g(P_y)=0$}}\right.\right\}.\]
Then $\dim_{\R} \cH_1= 2$ and $\cP_{3,3}^+ \cap \cH_1
  = \R_+ \cdot \frh_{p,0} + \R_+ \cdot x(x+p y-p z)^2$. 
Since $f$ is irreducible, $\cL_f = \R_+ \cdot \frh_{p,0}$. 

\smallskip

(4) Consider the case ($Cy$) and ($Cz$) are false. 

$f(P_y)=f(P_z)=0$ is impossible by \Ttagb. 
We may assume $f(P_y) > 0$. 
Then $G(x,y,z) = y(y-z)^2 \in \cL_f = \R_+ \cdot f$ as (2). 
A contradiction. 
\end{proof}

\bigbreak
\subsection{The case $V_{\C}(f)$ has a node}\hfill\break\relax
\hbox{}\hskip\parindent 
Assume that $f \in \cE(\cP_{3,3}^+)$ is irreducible 
and $V_{\C}(f)$ does not have a cusp. 
Then $V_{\C}(f)$ has the unique node $P$ on $\P_+^2$ by \Ttaea. 
$P$ must be an acnode. 
By \Ttah, $P \in \Int(\P_+^2)$. 
Let $P = (a \colon b \colon 1)$ ($a>0$, $b>0$). 
Note that if we replace the coordinate system $(x \colon y \colon z)$ by 
$(x/a \colon y/b \colon z)$, we can assume $P=(1 \colon 1 \colon 1)$. 

\def\Ttakb{Lemma 3.14}
\proclaim{Lemma 3.14} 
{\sl Put $Q_x := (0 \colon 1 \colon p)$, $Q_y := (1 \colon 0 \colon q)$, 
$P := (1 \colon 1 \colon 1)$ and }
\begin{align*}
 & \frg_{pq}(x,y,z) := z^3 + q^2 x^2z + p^2 y^2z - 2q x z^2 - 2p y z^2 
             - (p^2+q^2-4p-4q+3)x y z \\
 & \hskip70pt + (1-p+q)(1-p-q)x^2y + (1+p-q)(1-p-q) x y^2.
\end{align*}
{\parindent=20pt
\Item{(1)} {\sl Assume that $p>0$ and $q>0$. 
If $f \in \cP_{3,3}^+$ satisfies 
\[f(P) = f(Q_x) = f(Q_y) = f(P_x) = f(P_y) = 0,\]
then there exists $\alpha \in \R$ such that $f = \alpha \frg_{pq}$. }
\Item{(2)} {\sl Assume that $p>0$ and $q=0$. 
If $f \in \cH_{3,3}$ satisfies 
\[f(P) = f(Q_x) = f(P_x) = f(P_y) = f_z(P_y) = f_{zz}(P_y) = 0,\]
then there exists $\alpha \in \R$ such that $f = \alpha \frg_{p0}$. }
\Item{(3)} {\sl Assume that $p=0$ and $q=0$. 
If $f \in \cH_{3,3}$ satisfies 
\[f(P) = f(P_x)=f_z(P_x)=f_{zz}(P_x) = f(P_y)=f_z(P_y)=f_{zz}(P_y) = 0,\]
then there exists $\alpha \in \R$ such that $f = \alpha \frg_{00}$. }

} 
\endproclaim

\Proof
(1) Since $f(1,1,1) = f(0,1,p) = f(1,0,q) = f(1,0,0) = f(0,1,0) = 0$ and 
$f \in \cH_{3,3}$, $f$ must satisfies 
$f(1,1,1)=0$, $f_x(1,1,1)=0$, $f_y(1,1,1)=0$, 
$f(0,1,p)=0$, $f_y(0,1,p)=0$, 
$f(1,0,q)=0$, $f_x(1,0,q)=0$, 
$f(1,0,0)=0$, and $f(0,1,0)=0$. 
Using Mathematica, we have that the solution space is equal to 
$\R \cdot \frg_{pq}(x,y,z)$. 
This implies 
\[\cL_P \cap \cL_{Q_x} \cap \cL_{Q_y} 
 \cap \cL_{P_x} \cap \cL_{P_y} \subset \R_+ \cdot \frg_{pq}.\]
If $\frg_{pq} \in \cP_{3,3}^+$, 
then $\cL_f = \R_+ \cdot \frg_{pq}$, 
and $\frg_{pq} \in \cE(\cP_{3,3}^+)$. 

We can prove (2) and (3) by the similay way as (1). 
\end{proof}

\def\Ttaka{Theorem 3.15}
\def\Ttakan{3.15}
\proclaim{Theorem 3.15} 
{\sl Assume that $p \geq 0$ and $q \geq 0$. \par
{\parindent=20pt
\Item{\rm (1)} If $p+q > 1$, 
then $\frg_{pq} \notin \cP_{3,3}^+$. 
\Item{\rm (2)} If $p+q = 1$, 
then $\frg_{pq}(x,y,z) = (p x + q y - z)^2 z \in \cE(\cP_{3,3}^+)$. 
\Item{\rm (3)} If $p+q < 1$, then $\frg_{pq} \in \cE(\cP_{3,3}^+)$. 
Moreover, $\frg_{pq}(x^2$, $y^2$, $z^2) \in \cE(\cP_{3,6})$. 

}}
\endproclaim

\Proof
(1) Assume that $p+q>1$. Let $a := 1/(2(p+q-1)) > 0$. 
Since 
\[\frg_{pq}(x,x,1) = (x-1)^2\big(1-2(p+q-1)x\big),\]
we have $\frg_{pq}(x,x,1)<0$ if $x>a$. 
Thus $\frg_{pq} \notin \cP_{3,3}^+$. 

\smallskip

(2) is easy to see. 

\smallskip

(3) We shall show that $\frg_{pq} \in \cE(\cP_{3,3}^+)$ 
if $p+q<1$. 

Let $g(x,y,z) := (1-p+q)x + (1+p-q)y - 2z$. 
Note that $g(1,1,1) = 0$. 
Fix $(x \colon y \colon z) \in \P_+^2$. 

Consider the case $g(x,y,z) \geq 0$. 
Then $x \geq z$ or $y \geq z$, since $g(1,1,1)=0$. 
If $x \geq z$, then 
\[\frg_{pq}(x,y,z) = (1-p-q)y(x-z)g(x,y,z) 
   + \big(q x + (1-q)y - z\big)^2 z \geq 0. \eqno (\Ttakan.1)\]
If $y \geq z$, then 
\[\frg_{pq}(x,y,z) = (1-p-q)x(y-z)g(x,y,z) 
   + \big((1-p)x + p y - z\big)^2 z \geq 0. \eqno (\Ttakan.2)\]
Consider the case $g(x,y,z) < 0$. 
Then $x \leq z$ or $y \leq z$, since $g(1,1,1)=0$. 
If $x \leq z$, then $\frg_{pq}(x,y,z) \geq 0$ by (\Ttakan.1). 
If $y \leq z$, then $\frg_{pq}(x,y,z) \geq 0$ by (\Ttakan.2). 
Thus $\frg_{pq}(x,y,z) \geq 0$ for 
all $(x \colon y \colon z) \in \P_+^2$, if $p+q<1$. 

It $p+q \leq 1$, $p \geq 0$ and $q \leq 0$, then 
$\frg_{pq} \in \cE(\cP_{3,3}^+)$ by \Ttakb. 

\smallskip

(3-2) We shall show that $\frg_{pq}$ is irreducible if $p+q<1$. 
Then $V_{\R}(f) - \{P\}$ has a parametrization 
\begin{align*}
  & x(t) := -\frac{(1-q+q t)^2}{t(1-p-q)((1+p-q)+t(1-p+q))}, \\
  & y(t) := -\frac{(p+t-p t)^2}{(1-p-q)((1+p-q)+t(1-p+q))}. 
\end{align*}
This implies $V_{\C}(f)$ is an irreducible rational cubic curve. 

\smallskip

(3-3) $\frg_{pq}(x^2$, $y^2$, $z^2) \in \cE(\cP_{3,6})$ follows 
from 
$N_P(\frg_{pq}) + N_{Q_x}(\frg_{pq}) + N_{Q_y}(\frg_{pq}) 
 + N_{P_x}(\frg_{pq}) + N_{P_y}(\frg_{pq}) 
 = 4+2+2+1+1 = 10$. 
\end{proof}

Note that 
\begin{align*}
 & \frg_{p0}(x,y,z) = (1-p)^2 x^2y + (1-p)^2 x y^2 + p^2 x^2z 
      - 2p x z^2 - (1-p)(3-p)x y z + z^3,\\
 & \frg_{00}(x,y,z) = x^2y + x y^2 + z^3 - 3x y z.
\end{align*}

\def\Ttala{Definition 3.16}
\proclaim{Definition 3.16} 
{\rm We define $\frf_{pqr}(x,y,z)$ as the following: }
\begin{align*}
 & a_1(p,q) := p q-p+1, \\
 & a_2(p,q,r) := p^2q r-p^2r+p q r-p q+2p r+p-r+1, \\
 & c_1(p,q,r) := q^2 a_1(p,q) a_1(r,p) a_2(p,q,r), \\
 & c_2(p,q,r) := - a_1(p,q)(2p^3q^3r^3 - 2p^3q^2r^3 + 6p^2q^2r^3 
       - 2p q^3r^3 + 3p q^3r^2 - 6p q^2r^3 \\
 & \hskip70pt + 3p q^2r^2 + 2q^2r^3 - p q^3 - 3q^2r^2 
       + 3p q^2 - 3p q + q^2 + p - 1), \\
 & c_3(p,q,r) := r a_1(p,q) (p^3q^3r^3 - p^3q^2r^3 + 3p^2q^2r^3 
       - p q^3r^3 - 3p q^2r^3 + 3p q^3r + q^2r^3 \\
 & \hskip70pt - 2p q^3 - 3p q^2r + 6p q^2 - 3q^2r - 6p q + 2q^2 + 2p - 2), \\
 & c_4(p,q,r):= -c_1(p,q,r)-c_1(q,r,p)-c_1(r,p,q)
        -c_2(p,q,r)-c_2(q,r,p) \\
 & \hskip70pt -c_2(r,p,q)-c_3(p,q,r)-c_3(q,r,p)-c_3(r,p,q), \\
 & \frf_{pqr}(x,y,z) := c_1(p,q,r)x^3 + c_1(q,r,p)y^3+ c_1(r,p,q)z^3 \\
 & \hskip70pt  + c_2(p,q,r)x^2y + c_3(p,q,r)x y^2 
               + c_2(q,r,p)y^2z + c_3(q,r,p)y z^2 \\
 & \hskip70pt  + c_2(r,p,q)z^2x + c_3(r,p,q)z x^2 + c_4(p,q,r)x y z. 
\end{align*}
\endproclaim

\def\Ttalb{Lemma 3.17}
\proclaim{Lemma 3.17} 
{\sl Put $Q_x := (0 \colon p \colon 1)$, $Q_y := (1 \colon 0 \colon q)$, 
$Q_z := (r \colon 1 \colon 0)$, $P := (1 \colon 1 \colon 1)$.} \par
{\parindent=20pt
\Item{(1)} {\sl Assume that $p>0$, $q>0$ and $r>0$. 
If $f \in \cP_{3,3}^+$ satisfies 
\[f(P) = f(Q_x) = f(Q_y) = f(Q_z) = 0,\]
then there exists $\alpha \in \R$ such that $f = \alpha \frf_{pqr}$. }
\Item{(2)} {\sl Assume that $p>0$, $q>0$ and $r=0$. 
If $f \in \cP_{3,3}^+$ satisfies 
\[f(P) = f(Q_x) = f(Q_y) = f(P_y) = f_x(P_y) = 0,\]
then there exists $\alpha \in \R$ such that $f = \alpha \frf_{p,q,0}$. }
\Item{(3)} {\sl Assume that $p>0$, $q=0$ and $r=0$. 
If $f \in \cP_{3,3}^+$ satisfies 
\[f(P) = f(Q_x) = f(P_x) = f_z(P_x) = f(P_y) = f_x(P_y) = 0,\]
then there exists $\alpha \in \R$ such that $f = \alpha \frf_{p,0,0}$. }
\Item{(4)} {\sl Assume that $p=0$, $q=0$ and $r=0$. 
If $f \in \cP_{3,3}^+$ satisfies 
\[f(P) = f(P_z) = f_y(P_z) = f(P_x) = f_z(P_x) = f(P_y) = f_x(P_y) = 0,\]
then there exists $\alpha \in \R$ such that $f = \alpha \frf_{0,0,0}$. }

}
\endproclaim

\Proof
$f \in \cP_{3,3}^+$ and $f(P)=0$ implies $f_x(P)=f_y(P)=0$. 
If $r>0$, $f(Q_z)=0$ implies $f(Q_z)=f_x(Q_z)=0$. 
So, in any case of (1), (2), (3), (4), $f$ must satisfy 
$f(1,1,1)=0$, $f_x(1,1,1)=0$, $f_y(1,1,1)=0$, 
$f(0,p,1)=0$, $f_y(0,p,1)=0$, 
$f(1,0,q)=0$, $f_z(1,0,q)=0$, 
$f(r,1,0)=0$, $f_x(r,1,0)=0$, 
where $p \geq 0$, $q \geq 0$ and $r \geq 0$. 
We know that the solution of the above equations is a multiple 
of $\frf_{pqr}(x,y,z)$, using Mathematica. 
\end{proof}

We shall study the conditions on $p$, $q$, $r$ 
for $\frf_{pqr} \in \cP_{3,3}^+$. 

\def\Ttalc{Lemma 3.18}
\proclaim{Lemma 3.18} 
(I) {\sl Assume that $p>0$, $q>0$ and $r>0$. Then, 
$\frf_{pqr}(1,0,0) > 0$, $\frf_{pqr}(0,1,0) > 0$, 
and $\frf_{pqr}(0,0,1) > 0$, 
if and only if $a_1(p,q)>0$, $a_1(q,r)>0$ and $a_1(r,p)>0$. }

(II) {\sl Assume that $p \geq 0$, $q \geq 0$ and $r \geq 0$. Then, 
$\frf_{pqr}(1,0,0) \geq 0$, $\frf_{pqr}(0,1,0) \geq 0$, 
and $\frf_{pqr}(0,0,1) \geq 0$, 
if and only if $a_1(p,q) \geq 0$, $a_1(q,r) \geq 0$ and $a_1(r,p) \geq 0$. }
\endproclaim

\Proof
(1) Assume that $a_1(p,q)>0$, $a_1(q,r)>0$, and $a_1(r,p)>0$. 

Then $a_2(p,q,r) = p r a_1(p,q) + p a_1(q,r) + a_1(r,p) > 0$. 
Thus, 
\[\frf_{pqr}(1,0,0) = c_1(p,q,r) = q^2 a_1(p,q) a_1(r,p) a_2(p,q,r) > 0.\]
Similarly, we have $\frf_{pqr}(0,1,0) > 0$ and 
$\frf_{pqr}(0,0,1) > 0$. 

\smallskip

(2) Assume that $\frf_{pqr}(1,0,0) > 0$, $\frf_{pqr}(0,1,0) > 0$, 
and $\frf_{pqr}(0,0,1) > 0$. 
We shall derive a contradiction assuming $a_1(p,q)<0$. 

(2-1) we shall show that $p>1$, $0<q<1$, $a_1(q,r) > 0$ and $a_1(r,p) > 0$. 

If $p \leq 1$, then $a_1(p,q) = p q+(1-p) \geq 0$. 
If $q \geq 1$, then $a_1(p,q) = p(q-1)+1 \geq 0$. 
Thus, if $a_1(p,q)<0$, then $p>1$ and $0<q<1$. 
Then, $a_1(q,r) = q r+(1-q)>0$, and $a_1(r,p) = r(p-1)+1>0$. 

\smallskip

Let 
\begin{align*}
 & b_1(p,q) := -p^2q + p^2 - pq - 2p + 1 = (p-1)^2 - p(p+1)q, \\
 & b_2(p,q) := p^2q^2 - 2 p^2 q - 2 p q + p^2 - 2p + 1, \\
 & r_0(p,q) := \frac{1 + p(1-q)}{b_1(p,q)}, \\
 & r_2(p,q) := \frac{p(1-q)^2-1-q}{q(pq+(p-1))}. 
\end{align*}
(2-2) We shall show that $b_1(p,q)>0$ and $r_0(p,q) < r < r_2(p,q)$. 
Since 
\begin{align*}
 & 0 < \frf_{pqr}(1,0,0) = q^2 a_1(p,q) a_1(r,p) a_2(p,q,r), \\
 & 0 < \frf_{pqr}(0,1,0) = r^2 a_1(q,r) a_1(p,q) a_2(q,r,p), \\
 & 0 < \frf_{pqr}(0,0,1) = p^2 a_1(r,p) a_1(q,r) a_2(r,p,q), 
\end{align*}
we have $a_2(p,q,r) < 0$, $a_2(q,r,p) < 0$ and $a_2(r,p,q) > 0$. 
Note that $a_2(p,q,0)=1+p(1-q)>0$. 
Since $0 > a_2(p,q,r) = - b_1(p,q) r + a_2(p,q,0)$, 
we have $b_1(p,q) > a_2(p,q,0)/r > 0$. 
Thus, $a_2(p,q,r)$ is monotonically decreasing on $r$. 
The equation $a_2(p,q,r)=0$ on $r$ has just one root $r = r_0(p,q)$. 
Since $a_2(p,q,r)<0$, we have $r > r_0(p,q)$. 

Since $0 > a_2(q,r,p) = r q(p q+(p-1))  + (-p q^2 + 2p q - p + q + 1)$ 
is monotonically increasing on $r$, we have $r < r_2(p,q)$. 
Since $r_0(p,q)<r<r_2(p,q)$, we have 
\[0 < r_2(p,q) - r_0(p,q) 
   = \frac{-a_1(p,q) b_2(p,q)}{q (p q+(p-1)) b_1(p,q)}.\]
Thus $b_2(p,q)>0$.

\smallskip

Remember that 
$0 < a_2(r,p,q) = (p-1)q r^2 - ((p-1)-(p+2)q)r + (1-q)$. 
This is a quadratic function on $r$ with $(p-1)q > 0$. 
Note that 
\begin{align*}
 & a_2\big(r_0(p,q),p,q\big) 
    = \frac{2 a_1(p,q) b_2(p,q)}{\big((p-1)^2-pq(p+1)\big)^2} < 0, \\
 & a_2\big(r_2(p,q),p,q\big)
    = \frac{2a_1(p,q) b_2(p,q)}{\big((p-1)+p q\big)^2} < 0. 
\end{align*}
Thus, if $r_0(p,q) \leq r \leq r_2(p,q)$, 
then $a_2(r,p,q) < 0$. A contradiction. 
Thus we have $a_1(p,q)>0$. 

Similarly, we obtain $a_1(q,r)>0$ and $a_1(r,p)>0$. 

\smallskip

We can obtain (II), if we consider limits $p \to +0$, $q \to +0$ 
or $r \to +0$ in (I). 
\end{proof}

\def\Ttamb{Lemma 3.19}
\proclaim{Lemma 3.19} 
{\sl $\frf_{pqr}$ has the following properties:} 
{\parindent=20pt
\Item{\rm (1)} $\frf_{pqr}(z,x,y) = \frf_{qrp}(x,y,z)$, 
$\frf_{pqr}(y,z,x) = \frf_{rpq}(x,y,z)$. 
\Item{\rm (2)} $\displaystyle 
 \frf_{\frac{1}{p}\frac{1}{q}\frac{1}{r}}(x,y,z) 
    = \frac{1}{p^4 q^4 r^4}\frf_{p,r,q}(x, z, y)$. 
\Item{\rm (3)} {\sl If $a_1(p,q) = 0$, i.e. 
if $q = (p-1)/p$, then 
$$\frf_{p,(p-1)/p,r}(x,y,z) 
 = \frac{(1+(p-1)r)^4 z \big((p-1)x + y - p z)^2}{p^2}.$$

}
\endproclaim

\Proof
These follow from direct calculations using Mathematica. 
\end{proof}

\def\Ttama{Theorem 3.20}
\def\Ttaman{3.20}
\proclaim{Theorem 3.20} 
{\sl Assume that $p \geq 0$, $q \geq 0$ and $r \geq 0$. Then:} \par
{\parindent=20pt
\Item{(1)} {\sl $\frf_{pqr} \in \cP_{3,3}^+$ 
if and only if $a_1(p,q) \geq 0$, $a_1(q,r) \geq 0$ and $a_1(r,p) \geq 0$. }
\Item{(2)} {\sl $\frf_{pqr}$ is irreducible 
if and only if $a_1(p,q) > 0$, $a_1(q,r) > 0$ and $a_1(r,p) > 0$. }
\Item{(3)} {\sl If $a_1(p,q) > 0$, $a_1(q,r) > 0$ and $a_1(r,p) > 0$, 
then $\frf_{pqr}(x^2$, $y^2$, $z^2) \in \cE(\cP_{3,6})$. 

}
\endproclaim

\Proof
(i) `Only if part' of (1): 
Assume that $\frf_{pqr} \in \cP_{3,3}^+$. 
Then $\frf_{pqr}(1,0,0) \geq 0$, $\frf_{pqr}(0,1,0) \geq 0$, 
$\frf_{pqr}(0,0,1) \geq 0$. 
Thus $a_1(p,q) \geq 0$, $a_1(q,r) \geq 0$, and $a_1(r,p) \geq 0$, 
by \Ttalc (II). 

(ii) `Only if part' of (2) follows from \Ttamb (3). 

\smallskip

(iii) `If part' of (1) and `if part' of (2): 

If $p=q=r=0$, then $\frf_{000}(x,y,z) = x^2y + y^2z + z^2x - 3xyz$. 
In this case, $\frf_{000} \in \cP_{3,3}^+$, and $\frf_{000}$ is irreducible. 
So, we may assume $p \geq 0$, $q \geq 0$ and $r > 0$, by \Ttamb (1). 

(iii-1) Assume that $a_1(p,q)>0$, $a_1(q,r)>0$, and $a_1(r,p)>0$. 

Let $\ell_t \subset \P_{\R}^2$ be the line defined by $y-z = t(x-z)$ 
where $t \in \R$. 
The intersection point of $V_{\R}(\frf_{pqr})$ and $\ell_t$ 
($\ne (1 \colon 1 \colon 1)$) is given by 
$P(t) := \big(x_{pqr}(t) : y_{pqr}(t) : z_{pqr}(t)\big)$, where 
\begin{align*}
 & x_{pqr}(t) := a_1(q,r) \big(t+(p-1)\big)^2 
      \Big(r^2 a_1(p,q) a_2(q,r,p) t \\
 & \hskip70pt - \big((p^2q^2r^2+1) a_1(q,r) + 2q r a_1(r,p) 
                   + 2p q r^2 a_1(p,q)\big)\Big), \\
 & y_{pqr}(t) := a_1(r,p) \big((1-q)t+q\big)^2 
     \Big( - \big((p^2q^2r^2+1)a_1(r,p) \\
 & \hskip70pt +  2p r a_1(p,q) + 2p^2q r a_1(q,r)\big)t 
        + a_1(p,q) a_2(p,q,r) \Big), \\
 & z_{pqr}(t) := (r t-1)^2 a_1(p,q) 
      \big(a_1(q,r) a_2(q,r,p) t + q^2 a_1(r,p) a_2(p,q,r) \big). 
\end{align*}
Note that this also implies that $V_{\C}(\frf_{pqr})$ is irreducible. 
Thus we obtain the `if part' of (2). 
We shall show $\frf_{pqr} \in \cE(\cP_{3,3}^+)$. 

\smallskip

(iii-1-1) Consider the case $t \leq 0$. 
Then $x_{pqr}(t) \leq 0$ and $y_{pqr}(t) \geq 0$. 
Thus $P(t) \notin \Int(\P_+^2)$. 

(iii-1-2) Consider the case $t > 0$. Let 
\begin{align*}
 & t_1 := \frac{(p^2q^2r^2+1) a_1(q,r) + 2q r a_1(r,p) 
            + 2p q r^2 a_1(p,q)}{r^2 a_1(p,q) a_2(q,r,p)}, \\
 & t_2 := \frac{a_1(p,q) a_2(p,q,r)}{(p^2q^2r^2+1)a_1(r,p) + 
           2p r a_1(p,q) + 2p^2q r a_1(q,r)}. 
\end{align*}
$z_{pqr}(t) \geq 0$ for all $t > 0$. 
If $t < t_1$ then $x_{pqr}(t) \leq 0$, 
and if $t \geq t_1$ then $x_{pqr}(t) \geq 0$. 
If $t < t_2$ then $y_{pqr}(t) \geq 0$, 
and if $t \geq t_2$ then $y_{pqr}(t) \leq 0$. 
Note that 
\begin{align*}
 & t_1 - t_2 \\
 & \hskip10pt =
    \frac{a_2(r,p,q)\big(a_1(r,p) + p r a_1(p,q) + p^2q r a_1(q,r)big)}{
      r^2 a_1(p,q) a_2(q,r,p)
        \Big(a_1(r,p) + 2p r a_1(p,q) + 2p^2q r a_1(q,r) 
          + p^2q^2r^2a_1(r,p)\Big)} \\
 & \hskip10pt > 0. 
\end{align*}
Thus, for every $t > 0$, at least one of $x_{pqr}(t)$, 
$y_{pqr}(t)$, $z_{pqr}(t)$ is non-negative, 
and at least one of $x_{pqr}(t)$, $y_{pqr}(t)$, $z_{pqr}(t)$ is non-positive. 
Thus $P(t) \notin \Int(\P_+^2)$. 
Therefore $\frf_{pqr} \in \cP_{3,3}^+$. 
By \Ttalb, we have $\frf_{pqr} \in \cE(\cP_{3,3}^+)$. 

\smallskip

(iii-2) Consider the cases that at least one of 
$a_1(p,q)$, $a_1(q,r)$, $a_1(r,p)$ is equal to zero. 

$\frf_{pqr}(x,y,z)$ is continuous with respect to $a_1(p,q)$, 
$a_1(q,r)$ and $a_1(r,p)$. 
Thus, the `if part' of (2) follows from the above result. 

\smallskip

(3) $\frf_{pqr}(x^2$, $y^2$, $z^2) \in \cE(\cP_{3,6})$ follows 
from 
$N_P(\frf_{pqr}) + N_{Q_x}(\frf_{pqr}) + N_{Q_y}(\frf_{pqr}) 
 + N_{Q_z}(\frf_{pqr}) \allowbreak = 4+2+2+2 = 10$, 
where $Q_x := (0 \colon p \colon 1)$, $Q_y := (1 \colon 0 \colon q)$ and 
$Q_z := (r \colon 1 \colon 0)$ with $p \geq 0$, $q \geq 0$ and $r \geq 0$. 
\end{proof}

\bigbreak
\subsection{Final classification of $\cE(\cP_{3,3}^+)$}\hfill\break\relax
\hbox{}\hskip\parindent 
In this subsection, we essentially complete the classification 
of $f \in \cE(\cP_{3,3}^+)$ in the case 
$f$ has a node at $(1 \colon 1 \colon 1)$. 

\def\Ttapa{Corollary 3.21}
\proclaim{Corollary 3.21} 
{\sl Let $p>0$, $q>0$ and $r>0$. 
Assume that $f \in \cE(\cP_{3,3}^+)$ is irreducible, and satisfies 
\[f(1,1,1)=f(0,p,1)=f(1,0,q)=f(r,1,0)=0.\]
Then there exists $\alpha>0$ such that $f = \alpha \frf_{pqr}$. 
Moreover, $p q-p+1>0$, $q r-q+1>0$ and $r p-r+1>0$ hold. }
\endproclaim

\Proof
This follows from \Ttalb (1) and \Ttaman (2). 
\end{proof}

\def\Ttapb{Theorem 3.22}
\proclaim{Theorem 3.22} 
{\sl Let $P = (1 \colon 1 \colon 1)$, 
$Q_x = (0 \colon p \colon 1)$ and $Q_y = (q \colon 0 \colon 1)$. 
Assume that $f \in \cE(\cP_{3,3}^+)$ is irreducible, 
and satisfies $f(P)=f(Q_x)=f(Q_y)=0$ and $V_+(f) \cap L_z = \emptyset$. 
Then, one of the following statements holds: } \par
{\parindent=20pt
\Item{(1)} {\sl $f = \alpha \frg_{\frac{1}{p},\frac{1}{q}}$ 
($\exists \alpha > 0$) where $\displaystyle \frac{1}{p} + \frac{1}{q} < 1$. 
In this case $f(P_x)=f(P_y)=0$. }
\Item{(2)} {\sl $f(x,y,z) = \alpha \frf_{\frac{1}{p},0,q}(x,z,y)$ 
($\exists \alpha > 0$) 
where $p>1$ and $\displaystyle q<\frac{p}{p-1}$. 
In this case $f(P_x)=0$. }
\Item{(3)} {\sl $f = \alpha \frf_{p,\frac{1}{q},0}$ 
($\exists \alpha > 0$) where $q>1$ and $\displaystyle p<\frac{q}{q-1}$. 
In this case $f(P_y)=0$. }

}
\endproclaim

\Proof
Assume that $f \in \cE(\cP_{3,3}^+)$ is irreducible, 
$f(P)=f(Q_x)=f(Q_y)=0$ and $V_+(f) \cap L_z = \emptyset$. 
In \Ttagc, we put $P_1 = P$, $P_2=Q_x$ and $P_3=Q_y$. 
Let 
\[\cP_1 := \left\{ g \in \cP_{3,3}^+ \; \left| \; 
 \vcenter{
   \hbox{$g(P)=g_x(P)=g_y(P)=0$,}
   \hbox{$g(Q_x)=g_y(Q_x)=g(Q_y)=g_z(Q_y)=0$}}\right.\right\}.\]
Note that $\cP_1 = \cL_{P_1} \cap \cL_{P_2} 
  \cap \cL_{P_3}$. 
A direct calculation using Mathematica shows that $\dim \cP_1 = 3$. 
Thus $\{P$, $Q_x$, $Q_y\} \subsetne V_+(f)$. 
Since $V_+(f) \cap \Int \P_+^2 = \{P\}$, $V_+(f) \cap L_x = \{Q_x\}$, 
$V_+(f) \cap L_y = \{Q_y\}$ and $V_+(f) \cap L_z = \emptyset$, 
we obtain $\{P_x$, $P_y$, $P_z\} \cap V_+(f) \ne \emptyset$. 
Then, there are following four cases (i)-(iv). 

\smallskip

(i) The case $f(P_z)=0$. 

Then $f(P_x)>0$ and $f(P_y) > 0$ by \Ttagb. 
Moreover $I_{P_z}(f$, $x)=1$ and $I_{P_z}(f$, $y)=1$. 
Thus $V_+(f) = \{P$, $Q_x$, $Q_y$, $P_z\}$. 
$\cP_2 := \cP_1 \cap (\cP_{3,3}^+)_{P_z}$ is a two dimensional 
fan whose edges are generated by extremal elements 
$y((p-1)x + y - p z)^2$ and $x(q x + (1-q)y - z)^2$. 
Thus $\cP_2 \not\supset \R_+ \cdot f$. 
Therefore, we have $f(P_z) > 0$. 

Note that $I_{P_x}(f$, $y) \leq 1$, $I_{P_y}(f$, $x) \leq 1$  by \Ttagb. 

\smallskip

(ii) The case $f(P_x)=f(P_y)=0$. 
Then, there exists $\alpha \in \R$ such that 
$f = \alpha \frg_{\frac{1}{p},\frac{1}{q}}$, by \Ttakb(1). 

\smallskip

(iii) The case $f(P_x)=0$ and $f(P_y)>0$. 
Let 
\[\cP_3 := \big\{ g \in \cP_1 \; \big| \; 
   \hbox{$g(P_x) = g_y(P_x) = 0$}\big\}.\]
Then $\cP_3 = \R_+ \cdot \frf_{\frac{1}{p},0,q}(x,z,y)$ by \Ttalb. 
Since $I_{P_x}(f$, $y) \leq 1$, we have $g_z(P_x) \ne 0$. 
Note that if $g_x(P_x) = 0$ then $g_y(P_x) = 0$ by \Tqbaib. 
This implies that $\R_+ \cdot f = \cL_f \supset \cP_3$. 
Therefore, $f(x,y,z) = \alpha \frf_{\frac{1}{p},0,q}(x,z,y)$. 

\smallskip

(iv) The case $f(P_x)>0$ and $f(P_y)=0$. 
Same as (ii). 
\end{proof}


Assume that $f \in \cE(\cP_{3,3}^+)$ is irreducible
and satisfies $f(0$, $p$, $1) = 0$ for $p>0$. 
Then, $f(P_y)>0$ or $f(P_z)>0$ by \Ttagb. 

\def\Ttapc{Theorem 3.23}
\proclaim{Theorem 3.23} 
{\sl Let $P = (1 \colon 1 \colon 1)$ and $Q_x = (0 \colon p \colon 1)$. 
Assume that $f \in \cE(\cP_{3,3}^+)$ is irreducible, 
and satisfies $f(P)=f(Q_x)=0$, $f(P_z)>0$, $V_+(f) \cap L_y = \emptyset$ 
and $V_+(f) \cap L_z = \emptyset$. 
Then, $f(P_x)=f(P_y)=0$ and one of the following statements holds: } \par
{\parindent=20pt
\Item{(1)} {\sl $f = \alpha \frg_{\frac{1}{p},0}$ 
($\exists \alpha > 0$). }
\Item{(2)} {\sl $f(x,y,z) = \alpha \frf_{p,0,0}(x,y,z)$ 
($\exists \alpha > 0$) where $p<1$. }

}
\endproclaim

\Proof
Assume that $f \in \cE(\cP_{3,3}^+)$ is irreducible, 
and $V_+(f) \subset \{P$, $Q_x$, $P_x$, $P_y\}$. 

In \Ttagc, we put $P_1 = P$ and $P_2=Q_x$. 
Let 
\[\cP_1 := \left\{ g \in \cP_{3,3}^+ \; \left| \; 
 \vcenter{
   \hbox{$g(P)=g_x(P)=g_y(P)=0$,}
   \hbox{$g(Q_x)=g_y(Q_x)=0$}}\right.\right\}.\]
Note that $\cP_1 = \cL_{P_1} \cap \cL_{P_2}$. 

For $u \ne v \in \{x$, $y$, $z\}$, we denote $I_{uv} := I_{P_u}(f$, $v)$. 
For example, $I_{xy} = I_{P_x}(f$, $y)$. 
Let $I_x := \max \{I_{xy}$, $I_{xz}\}$, $I_y := \max \{I_{yz}$, $I_{yx}\}$ and 
$I_z := \max \{I_{zx}$, $I_{zy}\}$. 
Note that $I_x = \length_{P_x} f$, 
$I_y = \length_{P_y} f$ and $I_z = \length_{P_z} f$. 

\smallskip

(i) Consider the case $I_x + I_y = 4$. 
Then $(I_x$, $I_y) = (3$, $1)$, $(2$, $2)$ or $(1$, $3)$. 

\smallskip

(i-1) The case $I_x = 3$ and $I_y = 1$. 

Since $I_{xz} + I_{yz} \leq \deg f = 3$ and $I_{yz} \ne 0$, 
we have $I_{xz} \leq 2$. 
Thus $I_{xy} = 3$ and $I_{xz}=1$. Let 
\[\cP_2 := \big\{g \in \cP_1 \; \big| \;
  \hbox{$f(P_x)=f_z(P_x)=f_{zz}(P_x)=f(P_y)=0$}\big\}.\]
Then $\cP_2 = \R_+ \cdot \frg_{1/p,0}$ by \Ttakb(2). 
Since $\R_+ \cdot f = \cL_f \subset \cP_2$, 
we have $f = \alpha \frg_{1/p,0}$ ($\exists \alpha>0$). 

\smallskip

(i-2) The case $I_x = I_y = 2$. 

Since $I_{yx}=1$, we have $I_{yz}=2$. 
Since $I_{yz}+I_{xz} \leq 3$, we have $I_{xz}=1$ and $I_{xy}=2$. Let 
\[\cP_3 := \big\{g \in \cP_1 \; \big| \;
  \hbox{$f(P_x)=f_z(P_x)=f(P_y)=f_x(P_y)=0$}\big\}.\]
Then $\cP_3 = \R_+ \cdot \frf_{p,0,0}$ by \Ttalb(3). 
Since $\R_+ \cdot f = \cL_f \subset \cP_3$, 
we have $f = \alpha \frf_{p,0,0}$ ($\exists \alpha>0$). 

\smallskip

(i-3) The case $I_x = 1$ and $I_y = 3$. 

Since $I_{xz}+I_{yz} \leq 3$ and $I_{xz}=I_{xy}=I_x=1$, 
we have $I_{yz}=1$ and $I_{xy}=3$. Let 
\[\cP_4 := \big\{g \in \cP_1 \; \big| \;
  \hbox{$f(P_x)=f(P_y)=f_x(P_y)=f_{xx}(P_y)=0$}\big\}.\]
Then $\cP_4 = \R_+ \cdot z((p-1)x + y - p z)^2$. 
Thus, $f$ is reducible. 

\smallskip

(ii) Consider the case $I_x + I_y \ne 4$. 

If $I_x + I_y \leq 3$, then 
$\cL_f \supset \cP_2$, $\cL_f \supset \cP_3$ 
or $\cL_f \supset \cP_4$. 
If $I_x + I_y \geq 5$, then 
$\cL_f \subset \cP_2$, $\cL_f \subset \cP_3$ 
or $\cL_f \subset \cP_4$. 
In any case, $\cL_f = \cP_2$ or $\cL_f = \cP_3$ 
or $\cL_f = \cP_4$. 
\end{proof}

\def\Ttapd{Theorem 3.24}
\proclaim{Theorem 3.24} 
{\sl Assume that $f \in \cE(\cP_{3,3}^+)$ is irreducible, 
and satisfies $f(P)=0$ where $P = (1 \colon 1 \colon 1)$. 
Moreover, we assume that $V_+(f) \subset \{P$, $P_x$, $P_y$, $P_z\}$ 
and $\length_{P_x} f \geq \length_{P_y} f \geq \length_{P_z} f$. 
Then one of the following statements holds: } \par
{\parindent=20pt
\Item{(1)} {\sl $f = \alpha \frg_{0,0}$ ($\exists \alpha > 0$).} 
\Item{(2)} {\sl $f = \alpha \frf_{0,0,0} = \alpha(x^2y+y^2z+z^2x-3xyz)$ 
($\exists \alpha > 0$). }
\Item{(3)} {\sl $f(x,y,z) = \alpha \frf_{0,0,0}(x,z,y) 
 = \alpha(xy^2+yz^2+zx^2-3xyz)$ ($\exists \alpha > 0$). }

}
\endproclaim

\Proof
Assume that $f \in \cE(\cP_{3,3}^+)$ is irreducible, 
and $V_+(f) \subset \{P$, $P_x$, $P_y$, $P_z\}$. 
In \Ttagc, we assume $P_1 = P$. 
Let $I_x$, $I_y$, $I_z$, $I_{xy}$,$\ldots$, $I_{zy}$ be same as in 
the proof of \Ttapc. 
By our assumption, $I_x \geq I_y \geq I_z$. 
If $I_{xy} \geq 2$, then $I_{xz} = 1$. 
Note that $I_{yx}+I_{zx} \leq 3$, $I_{xy}+I_{zy} \leq 3$ and 
$I_{xz}+I_{yz} \leq 3$ by \Ttagb. 
In particular, $I_x+I_y+I_z \leq 9-1-1-1 = 6$. 

\smallskip

(i) Consider the case $I_x + I_y + I_z = 6$. 
Then $(I_x$, $I_y$, $I_z) = (3$, $3$, $0)$, 
$(3$, $2$, $1)$ or $(2$, $2$, $2)$. 
By \Ttagb, $(I_x$, $I_y$, $I_z) = (3$, $2$, $1)$ is impossible. 

\smallskip

(i-1) The case $I_x = I_y = 3$ and $I_z = 0$. 

This can possible only if $I_{xy} = I_{yx} = 3$. Thus $f$ is an element of 
\[\cP_1 := \left\{ g \in \cP_{3,3}^+ \ \left| \ 
  \vcenter{ 
     \hbox{$g(P)=g(P_x)=g_z(P_x)=g_{zz}(P_x)=0$,}
     \hbox{$g(P_y)=g_z(P_y)=g_{zz}(P_y)=0$}}\right.\right\}.\]
$\cP_1 = \R_+ \cdot \frg_{00}$ \Ttakb(3). 
Thus $f = \alpha \frg_{00}$ ($\exists \alpha > 0$). 

\smallskip

(i-2) The case $I_x = I_y = I_z = 2$. 

By \Ttakb(3), there are two possibilities that 
$I_{xy} = I_{yz} = I_{zx} = 2$ and $I_{xz} = I_{yx} = I_{zy} = 2$. 

If $I_{xy} = I_{yz} = I_{zx} = 2$, then $f$ is an element of 
\[\cP_2 := \left\{ g \in \cP_{3,3}^+ \ \left| \ 
  \vcenter{
      \hbox{$g(P)=g(P_x)=g_z(P_x)=0$,}
      \hbox{$g(P_y)=g_x(P_y)=g(P_z)=g_y(P_z)=0$}}\right.\right\}.\]
Then, $\cP_2 = \R_+ \cdot \frf_{000}$ by \Ttalb(4), 
and we have (2). 

If $I_{xz} = I_{yx} = I_{zy} = 2$, then $f$ is an element of 
\[\cP_3 := \left\{ g \in \cP_{3,3}^+ \ \left| \ 
  \vcenter{
    \hbox{$g(P)=g(P_x)=g_y(P_x)=0$,}
    \hbox{$g(P_y)=g_z(P_y)=g(P_z)=g_x(P_z)=0$}}\right.\right\}.\]
Then, $\cP_3 = \R_+ \cdot \frf_{000}(x,z,y)$, and we have (3). 

\smallskip

(ii) Consider the case $I_x + I_y + I_z \leq 5$ and $I_{xy} > I_{xz}$. 

Then $I_y \leq 2$, $I_z \leq 1$ and $I_{xz}=1$. 
Let $G(x,y,z) = y(x-z)^2$. 
Then $G_x=2y(x-z)$, $G_y(x-z)^2$, $G_z=-2y(x-z)$, $G_{xx}=G_{zz}=2y$ 
and $G_{yy} = 0$. Note that $G_y(P_x) \ne 0$ and $G_y(P_z) \ne 0$. 

Since $G(P)=G_x(P)=G_y(P)=G_z(P)=0$, we have $G \in \cL_P$. 
Since $G(P_x)=G_z(P_x)=G_{zz}(P_x)=0$, we have $G \in \cL_{P_x}$. 
Since $G(P_y)=G_x(P_y)=G_z(P_y)=0$, we have $G \in \cL_{P_y}$. 
Since $G(P_z)=0$, we have $G \in \cL_{P_z}$. 
Thus $G \in \cL_f = \R_+ \cdot f$. 
A contradiction. 

\smallskip

(iii) Consider the case $I_x + I_y + I_z \leq 5$ and $I_{xy} < I_{xz}$. 

As the above argument, we have 
$G(x,z,y) = z(x-y)^2 \in \cL_f = \R_+ \cdot f$. 
A contradiction. 

\smallskip

(iv) Consider the case $I_x + I_y + I_z \leq 5$ and $I_{xy} = I_{xz}$. 

Then $I_x=I_y=I_z=1$, and Thus $G \in \cL_f = \R_+ \cdot f$. 

Thus, we complete the proof. 
\end{proof}

\bigbreak
\subsection{Proofs of Theorems in \S 1}\hfill\break\relax
\hbox{}\hskip\parindent 
\Ttac \ follows from \Ttama (2) and \Ttalb (1). 
\Ttad \ follows from \Ttaka (2), (3) and \Ttakb (1). 
\Ttadb \ follows frrom \Ttcb \ and \Ttca. 

\bigskip

{\it Proof of \Ttae. }
(I) If $f \in \cE(\cP_{3,3}^+)$ is reducible, 
then $f$ is (4) or (5) in \Ttae, by \Ttaca.

Assume that $f \in \cE(\cP_{3,3}^+)$ is irreducible, 
Then, $V_{\C}(f)$ is a rational curve on $\P_{\C}^2$ whose unique 
singular point lies on $\P_+^2$. 

\smallskip

(II) Consider the case that $V_{\C}(f)$ has a cusp $P$. 
Then $P \in L_x \cup L_y \cup L_z$ by \Ttcb. 
If $P = (0 \colon a \colon 1) \in L_x$ ($a>0)$, 
put $f_1(x,y,z) = f(x,ay,z)$. 
Then, $V_{\C}(f_1)$ is a rational curve 
whose cusp is at $(0 \colon 1 \colon 1)$. 
Thus $f_1 = \alpha \frh_{pq}$ ($\exists \alpha>0$)by \Ttadb. 
So, $f(x$, $y$, $z) = \alpha \frh_{pq}(x$, $y/a$, $z)$. 

If $P = (1 \colon 0 \colon a) \in L_y$ ($a>0)$, 
put $f_1(x$, $y$, $z) = f(z$, $x$, $ay)$. 
Then, $V_{\C}(f_1)$ is a rational curve 
whose cusp is at $(0 \colon 1 \colon 1)$. 
Thus $f_1(x$, $y$, $z) = \alpha \frh_{pq}(x$, $y$, $z)$ and 
$f(x$, $y$, $z) = \alpha \frh_{pq}(y$, $z/a$, $x)$. 

Similarly, if $P = (a \colon 1 \colon 0) \in L_z$ ($a>0)$, 
then $f(x$, $y$, $z) = \alpha \frh_{pq}(z$, $x/a$, $y)$. 

\smallskip

(III) We consider the case that $V_{\C}(f)$ has an acnode $P$. 
Then $P \in \Int(\P_+^2)$ by \Ttah. 
Let $P = (1 \colon a \colon b)$ ($a>0$, $b>0$), and 
$f_1(x$, $y$, $z) := f(x$, $ay$, $bz)$. 
Then $f_1 \in \cE(\cP_{3,3}^+)$ and 
$V_{\C}(f_1)$ has an acnode at $(1 \colon 1 \colon 1)$. 
After a suitable permutation $\sigma$ of $x$, $y$, $z$, 
$f_1(\sigma(x)$, $\sigma(y)$, $\sigma(z))$ is equal to one of 
$\alpha \frf_{pqr}(x$, $y$, $z)$, 
$\alpha \frg_{\frac{1}{p},\frac{1}{q}}(x$, $y$, $z)$, 
$\alpha \frf_{\frac{1}{p},0,q}(x$, $z$, $y)$, 
$\alpha \frf_{p,\frac{1}{q},0}(x$, $y$, $z)$, 
$\alpha \frg_{\frac{1}{p},0}(x$, $y$, $z)$, 
$\alpha \frf_{p,0,0}(x$, $y$, $z)$, 
$\alpha \frg_{0,0}(x$, $y$, $z)$, 
$\alpha \frf_{0,0,0}(x$, $y$, $z)$ and 
$\alpha \frf_{0,0,0}(x$, $z$, $y)$ 
by results in \S 3.5. 
On the other hand 
\begin{align*} 
 & \frg_{pq}(y,x,z) = \frg_{qp}(x,y,z) \\
 & \frf_{pqr}(y,x,z) = \frf_{qrp}(x,z,y), \\ 
 & \frf_{pqr}(y,z,x) = \frf_{rpq}(x,y,z), \\
 & \frf_{pqr}(z,x,y) = \frf_{qrp}(x,y,z), \\
 & \frf_{pqr}(z,y,x) = \frf_{rpq}(x,z,y), \\
 & \frf_{pqr}(y,x,z) = p^4 q^4 r^4 
     \frf_{\frac{1}{q},\frac{1}{p},\frac{1}{r}}(x, y, z), 
\end{align*}
Thus $f(x$, $y$, $z)$ can be represented as one of the forms 
in (1), (2) and (3). \QED

\def\Ttanc{Corollary 3.25}
\proclaim{Corollary 3.25} 
{\rm (1)} $\cE(\cP_{3,6}^e) \subset \cE(\cP_{3,6})$. \par
{\parindent=20pt
\Item{(2)} {\rm (\Ttab)} {\sl If $f(x,y,z) \in \cE(\cP_{3,3}^+)$, 
then $f(x^2,y^2,z^2) \in \cE(\cP_{3,6})$. }

}
\endproclaim

\Proof
(1) and (2) are equivalent, 
since $\cP_{3,3}^+ \cong \cP_{3,6}^e$ by the correspondence 
$f(x,y,z) \allowbreak \to f(x^2,y^2,z^2)$. 

Take $g \in \cE(\cP_{3,6}^e)$. 
There exists a $f \in \cE(\cP_{3,3}^+)$ such that 
$g(x,y,z) = f(x^2$, $y^2$, $z^2)$. 
If $f(x,y,z) = \frf_{pqr}(x$, $y/a$, $z/b)$, 
then $g \in \cE(\cP_{3,6})$ by \Ttama (1). 
Similarly, we obtain $g \in \cE(\cP_{3,6})$ in the cases 
(2), (3), (4) of \Ttae, by \Ttca \ and \Ttaka. 

Consider the case $f(x,y,z)=x(ax+bx+cx)^2$ 
with $\dim \big(V_{\R}(ax+by+cz) \cap \P_+^2\big) = 1$. 
Then $V_{\R}(ax^2+by^+cz^2)$ is an irreducible real quadric curve. 
Thus $g(x,y,z)=x^2(ax^2+bx^2+cx^2)^2 \in \cE(\cP_{3,6})$. 

It is easy to see that $x^2y^2z^2 \in \cE(\cP_{3,6})$. 
Thus we have the conclusion. \QED
\end{proof}

\def\Ttana{Proposition 3.26}
\proclaim{Proposition 3.26} 
{\sl Let 
$\cP_{3,3}^{c+} := \big\{ f \in \cP_{3,3}^+$ $\big|$ 
$f(y,z,x) = f(x,y,z)$ (i.e. $f$ is cyclic) $\big\}$. 
Then, 
$\cE(\cP_{3,3}^{c+}) \subset \cE(\cP_{3,3}^+)$. }
\endproclaim

\Proof
By Theorem 3.1 of \cite{RefAa}, $f \in \cE(\cP_{3,3}^{c+})$ is equal to 
$cxyz$ or $c\frf_s$ ($\exists c > 0$, $s \in [0$, $+\infty]$) 
where 
\begin{align*}
 & \frf_s(x,y,z) := s^2 (x^3+x^3+z^3) - (2s^3-1)(x^2y+y^2z+z^2x) \\
 & \hskip70pt + (s^4-2s)(x y^2+y z^2+z x^2) - 3(s^4-2s^3+s^2-2s+1)x y z, \\
 & \frf_{\infty}(x,y,z) := x y^2+y z^2+z x^2 - 3x y z. 
\end{align*}
Note that $\frf_{sss}(x,y,z) = (s^2+s+1)(s^2-s+1)^2 \frf_s(x,y,z)$. 
Thus $\frf_s(x,y,z) \in \cE(\cP_{3,3}^+)$ for $s \geq 0$. 
Since $xyz \in \cE(\cP_{3,3}^+)$, 
we have $\cE(\cP_{3,3}^{c+}) \subset \cE(\cP_{3,3}^+)$. 
\end{proof}

\def\Ttanb{Corollary 3.27}
\proclaim{Corollary 3.27} 
{\sl All the cyclic elements of $\cE(\cP_{3,3}^{c+})$ are 
$c \frf_s$ ($s \in [0$, $+\infty]$) and $c xyz$ ($c > 0$). 
All the symmetric elements of $\cE(\cP_{3,3}^{c+})$ are 
$c \frf_1$ and $c xyz$ ($c>0$).}
\endproclaim

Note that $\frf_{t^2}(x^2$, $y^2$, $z^2)$ appeared 
as $S_t(x$, $y$, $z)$ in \cite[(1.9), (6.20)]{RefR}. 

\section{Appendix}
Hilbert proved $P_{3,4} = \Sigma_{3,4}$. 
We shall give an alternative proof of this theorem. 
The proof of \cite[Proposition 6.3.4]{RefBCR} is one of classical type 
proofs of this Hilbelt's theorem. 

\def\Ttda{Theorem 4.1}
\proclaim{Theorem 4.1}{\rm (Hilbert)} 
{\sl If $f \in \cE(\cP_{3,4})$, 
then $f$ is the square of a quadratic polynomial. }
\endproclaim

\Proof
Assume that $f \in \cE(\cP_{3,4})$ 
is not a square of a quadratic polynomial. 
It is easy to see that this implies that $f$ is irreducible in $\R[x,y,z]$. 

If $f$ is not exposed, then there exists $f_n \in \cE(\cP_{3,4})$ 
($n \in \N$) such that $\displaystyle \lim_{n \to \infty} f_n = f$ 
(with respect to the Euclidean topology of $\cH_{3,4}$), 
and that all $f_n$ are exposed (see \cite[Theorem 2.1.7]{RefSch}). 
If $f$ is irreducible in $\R[x,y,z]$, 
we can take $f_n$ so that $f_n$ are irreducible in $\R[x,y,z]$. 
So, we may assume that $f$ is exposed. 

Consider $V_{\R}(f)$. 
It is easy to see that if $\dim_{\R} V_{\R}(f) = 1$ as a topological space, 
then $f$ is not irreducible in $\R[x,y,z]$, since $f$ is PSD. 
If $V_{\R}(f) = \emptyset$, then $f$ cannot be extremal. 
Thus $V_{\R}(f)$ is a set of isolated points. 
Since $f$ is exposed, $V_{\R}(f)$ does not contain infinitely near points. 
Since $\dim_{\R} \cH_{3,4} = 15$, 
$V_{\R}(f)$ must contain at least $5 = 15/3$ points. 

If $V_{\C}(f)$ is irreducible, this is impossible 
because any curves on $\P_{\C}^2$ whose arithmetic genus is equal to 2, 
can have at most 4 singular points. 

Thus $f = g \overline{g}$ where $g \in \C[x,y,z] - \R[x,y,z]$ is a quadratic. 
Note that $V_{\R}(f) \subset V_{\C}(g) \cap V_{\C}(\overline{g})$. 
Thus $\# V_{\R}(f) \leq 4$. A contradiction. 
\end{proof}

If $f \in \cE(\cP_{3,6}) - \Sigma_{3,6}$, 
then $x^{2d}f \in \cE(\cP_{3,6+2d}) - \Sigma_{3,6+2d}$. 
Here we shall give examples of 
irreducible $f \in \cE(\cP_{3,d}) - \Sigma_{3,d}$ for $d=8$ and $10$. 
Such $f$ will be more interesting than reducible ones. 

\proclaim{Theorem 4.2} 
{\sl There exists $f(x,y,z) \in \cE(\cP_{3,4}^+)$ such 
that $f(x^2$, $y^2$, $z^2) \in \cE(\cP_{3,8})$ and 
$f(x^2$, $y^2$, $z^2)$ is irreducible in $\C[x,y,z]$. }
\endproclaim

\Proof
In $\cE(\cP_{3,4}^+)$, the equality conditions 
$f(1,1,1) = f(2,3,1) = f(1,2,3) = f(0,4,3) \allowbreak = f(6,0,5) \allowbreak 
= f(0,1,0) = 0$ determine the polynomial 
\begin{align*}
 f(x,y,z) & := 591900050 x^4 + 437205100 x^3 y - 766414561 x^2 y^2 
     + 217365672 x y^3 \\
  & \hskip15pt - 1650610670 x^3 z - 102695021 x^2 y z 
     + 248518503 x y^2 z + 549666 y^3 z \cr
  & \hskip15pt + 1531736792 x^2 z^2 + 118221267 x y z^2 
     + 101630538 y^2 z^2 \\
  & \hskip15pt - 636743352 x z^3 - 273946320 y z^3 
     + 183282336 z^4 
\end{align*}
up to a constant multiplication. 
We can prove that $f \in \cP_{3,4}^+$ if we observe 
$f(1+x$, $1+tx$, $1)/x^2$ carefully. 
It is easy to prove that $f(x^2$, $y^2$, $z^2)$ is irreducible 
in $\C[x,y,z]$ (see Theorem 4.3 below). 
$f(x^2$, $y^2$, $z^2)$ has the following $17$ isolated zeros: 
$(1 \colon 1 \colon 1)$, $(1 \colon 1 \colon -1)$, $(1 \colon -1 \colon 1)$, 
$(-1 \colon 1 \colon 1)$, 
$(\sqrt{2} \colon \sqrt{3} \colon 1)$, $(\sqrt{2} \colon \sqrt{3} \colon -1)$, 
$(\sqrt{2} \colon -\sqrt{3} \colon 1)$, 
$(-\sqrt{2} \colon \sqrt{3} \colon 1)$, 
$(1 \colon \sqrt{2} \colon \sqrt{3})$, $(1 \colon \sqrt{2} \colon -\sqrt{3})$, 
$(1 \colon -\sqrt{2} \colon \sqrt{3})$, 
$(-1 \colon \sqrt{2} \colon \sqrt{3})$, 
$(0 \colon 2 \colon \sqrt{3})$, $(0 \colon 2 \colon -\sqrt{3})$, 
$(\sqrt{6} \colon 0 \colon \sqrt{5})$, $(\sqrt{6} \colon 0 \colon -\sqrt{5})$ 
and $(0 \colon 1 \colon 0)$. 
Solve the differential equations for $F \in \cH_{3,8}$ such that 
$F(P) = F_x(P) = F_y(P) = 0$ for the above 17 points $P$. 
The solution space of this equation is $\R \cdot f(x^2,y^2,z^2)$. 
Thus $f(x^2$, $y^2$, $z^2) \in \cE(\cP_{3,8})$. 
\end{proof}

\proclaim{Theorem 4.3} 
{\sl There exists $f(x,y,z) \in \cE(\cP_{3,5}^+)$ such that 
$f(x^2$, $y^2$, $z^2) \in \cE(\cP_{3,10})$ 
and $f(x^2$, $y^2$, $z^2)$ is irreducible in $\C[x,y,z]$.}
\endproclaim

\Proof
The equality conditions 
$f(4,1,1) = f(1,4,1) = f(1,1,4) = f(1,9,9) = f(9,1,9) \allowbreak
 = f(9,9,1) = f(1,0,0) = f(0,1,0) = 0$ determine the polynomial 
\begin{align*}
 & f(x,y,z) = 837x^4y - 645x^3y^2 - 645x^2y^3 + 837x y^4 + 1755x^4z \\
 & \hskip60pt - 17181x^3y z + 23876x^2y^2z - 17181x y^3z 
     + 1755 y^4z - 3486x^3z^2 \\
 & \hskip60pt + 19594x^2y z^2 + 19594x y^2z^2 - 3486y^3z^2 + 3287x^2z^3 
     - 11030x y z^3 \\
 & \hskip60pt + 3287y^2z^3 - 1692x z^4 -1692y z^4 + 648z^5. 
\end{align*}
Elementary but somewhat long calculation shows that $f \in \cP_{3,5}^+$. 
Let $g(x$, $y$, $z) := f(x^2$, $y^2$, $z^2)$. 
Then $V_{\R}(g)$ has the following 26 acnodes: 
$V_{26} := \big\{(2 \colon 1 \colon 1)$, $(2 \colon 1 \colon -1)$, 
$(2 \colon -1 \colon 1)$, $(-2 \colon 1 \colon 1)$, 
$(1 \colon 2 \colon 1)$, $(1 \colon 2 \colon -1)$, 
$(1 \colon -2 \colon 1)$, $(-1 \colon 2 \colon 1)$, 
$(1 \colon 1 \colon 2)$, $(1 \colon 1 \colon -2)$, 
$(1 \colon -1 \colon 2)$, $(-1 \colon 1 \colon 2)$, 
$(1 \colon 3 \colon 3)$, $(1 \colon 3 \colon -3)$, 
$(1 \colon -3 \colon 3)$, $(-1 \colon 3 \colon 3)$, 
$(3 \colon 1 \colon 3)$, $(3 \colon 1 \colon -3)$, 
$(3 \colon -1 \colon 3)$, $(-3 \colon 1 \colon 3)$, 
$(3 \colon 3 \colon 1)$, $(3 \colon 3 \colon -1)$, 
$(3 \colon -3 \colon 1)$, $(-3 \colon 3 \colon 1)$, 
$(1 \colon 0 \colon 0)$, $(0 \colon 1 \colon 0) \big\}$. 
The solution space of the system of equalities 
$g(P)=g_x(P)=g_y(P)=g_z(P)=0$ for all $P \in V_{26}$, 
is equal to $\R \cdot f(x^2,y^2,z^2)$. 
Thus $g \in \cE(\cP_{3,10})$. 
So, $f \in \cE(\cP_{3,5}^+)$. 
We shall show that $g$ is irreducible in $\C[x,y,z]$. 

\smallskip

(1) To begin with, we prove that $g$ is irreducible in $\R[x,y,z]$. 
Assume that $g = h_1 h_2$ where $h_1 \in \cH_{3,d}$, 
$h_2 \in \cH_{3,e}$ with $d+e=10$, $d \leq e$. 
Since $g \in \cE(\cP_{3,10})$, 
we have $h_1 \in \cE(\cP_{3,d})$ and 
$h_2 \in \cE(\cP_{3,e})$. 
If $d$ is odd, then $\cP_{3,d} = 0$. 
Thus $d$ is even. 
If $d=2$, then $h_1 = h_3^2$ ($\exists h_3 \in \cH_{3,1}$). 
Then $V_{26}$ must contain a line $V_{\R}(h_3)$. 
If $d=4$, then $h_1 = h_4^2$ ($\exists h_4 \in \cH_{3,2}$). 
Since $\Sing(V_{\R}(h_4)) \subset V_{26}$, we have $V_{\R}(h_4) = \emptyset$ 
or $h_4 = h_5^2+h_6^2$ ($\exists h_5$, $h_6 \in \cH_{3,1}$). 
Is is easy to see that these are impossible. 

\smallskip 

(2) Thus, if $g$ is reducible, there exists an imaginary 
$h_7 \in \C[x,y,z]$ such that $g = h_7 \overline{h_7}$ 
where $\overline{h_7}$ is the complex conjugate of $h_7$. 
If $P \in \Sing(V_{\C}(h_7)) \cap \P_{\R}^2 \ne \emptyset$, 
then $P \in \Sing(V_{\C}(\overline{h_7})) \cap \P_{\R}^2$. 
This is impossible, since $P \in V_{26}$ is an acnode. 
Thus $\Sing(V_{\C}(h_7)) \cap \P_{\R}^2 = \emptyset$. 
This implies $V_{26} \subset V_{\C}(h_7) \cap V_{\C}(\overline{h_7})$. 
But $\# \big(V_{\C}(h_7) \cap V_{\C}(\overline{h_7})\big) \leq 5^2
  = 25$. 
\end{proof}


\end{document}